\documentclass[11pt, a4paper]{amsart}

\usepackage{amsfonts,amsmath,amssymb, amscd}

\usepackage[all]{xy}

\newtheorem{theorem}{Theorem}[section]
\newtheorem{lemma}[theorem]{Lemma}
\newtheorem{prop}[theorem]{Proposition}
\newtheorem{corollary}[theorem]{Corollary}

\theoremstyle{definition}
\newtheorem{definition}[theorem]{Definition}
\newtheorem{rem}[theorem]{Remark}

\newtheorem{example}[theorem]{Example}

\usepackage{mathrsfs}
\usepackage{color}

\newcommand\pf{\begin{proof}}
\newcommand\epf{\end{proof}}

\newcommand{\A}{A^{\natural}}
\newcommand{\bA}{\mathbb{A}}
\newcommand{\M}{\mathtt{SMod}}
\newcommand{\g}{\mathfrak{g}}
\newcommand{\Lie}{\operatorname{Lie}}

\newcommand{\ot}{\otimes}
\newcommand{\op}{\operatorname}
\newcommand{\grM}{\mathtt{GrMod}}
\newcommand{\cA}{\mathcal{A}}
\newcommand{\cH}{\mathcal{H}}
\newcommand{\cR}{\mathcal{R}}
\newcommand{\G}{\mathsf{G}}
\newcommand{\X}{\mathsf{X}}
\newcommand{\Y}{\mathsf{Y}}

\newcommand{\rcosmash}{\mathop{\raisebox{0.2ex}{\makebox[0.92em][l]{${\scriptstyle\blacktriangleright\mathrel{\mkern-4mu}<}$}}}}
\newcommand{\cO}{\mathcal{O}}
\newcommand{\cF}{\mathcal{F}}

\newcommand{\lboson}{\mathop{\raisebox{0.2ex}{\makebox[1.05em][l]{${\scriptstyle>\mathrel{\mkern-4mu}\lessdot}$}}\raisebox{0.12ex}{\hspace{-0.8mm}$\shortmid$}}}

\numberwithin{equation}{section}

\title[Torsors in super-symmetry]
{Torsors in super-symmetry}

\author[A.~Masuoka]{Akira Masuoka}
\address{Akira Masuoka,
Department of Mathematics, 
University of Tsukuba, 
Ibaraki 305-8571, Japan}
\email{akira@math.tsukuba.ac.jp}

\author[T.~Oe]{Takuya Oe}
\address{Takuya Oe,
Graduate School of Pure and Applied Sciences, 
University of Tsukuba, Ibaraki 305-8571, Japan;\
Present address:
NanoBridge Semiconductor, Inc.
Onogawa Technology Center, Onogawa, Tsukuba, Ibaraki 305-8569, Japan}

\author[Y.~Takahashi]{Yuta Takahashi}
\address{Yuta Takahashi,
Graduate School of Pure and Applied Sciences, 
University of Tsukuba, Ibaraki 305-8571, Japan;\
Present address:
National Institute of Technology, Numazu College,
3600 Ooka Numazu, Shizuoka 410-8501, Japan}
\email{y.takahashi@denki.numazu-ct.ac.jp}

\begin{document}

\begin{abstract}
Torsors under affine groups are generalized in the super context by 
super-torsors under affine super-groups. We investigate those super-torsors by using 
Hopf-algebra language and techniques.
It is shown in an explicit way that under suitable assumptions, every super-torsor arises from an 
ordinary torsor.
Especially, the objects with affinity restriction, 
or namely, 
the affine super-torsors and the affine ordinary torsors
are proved to be precisely in one-to-one correspondence. 
The results play substantial roles 
in ongoing construction of super-symmetric Picard-Vessiot theory.
\end{abstract}

\maketitle




\bigskip

\section{Introduction---basic notations and main theorems}\label{SINT}

Throughout in this paper we work over a fixed field $k$ whose characteristic 
$\operatorname{char}k$ differs from $2$, or in notation, $\operatorname{char}k\ne 2$. 

This section, which is slightly long as an introduction of a paper, is devoted mainly to
clarifying the notion of super-torsors and to formulating
two our main theorems on super-torsors. 

\subsection{Very basic algebraic super-geometry}\label{SSMiP}
Such basics are here presented mostly
by using the functor-of-points approach.

Recall that a $k$-\emph{functor} is a set-valued functor defined on the category 
$\mathtt{Alg}_k$ of commutative algebras (over $k$). A representable $k$-functor is called 
an \emph{affine scheme}, and it is called an \emph{affine group} if it is group-valued. 
In addition, the category of $k$-functors includes faisceaux (or sheaves), 
faisceaux dur (or dur sheaves),
schemes and affine schemes, as subcategories. These notions as well as relevant
basic notions and constructions
are directly generalized to the super context, 
in which everything is based on the
tensor category $\M_k$ of vector spaces graded by the order-2-group $\mathbb{Z}_2$,
equipped with the so-called \emph{super-symmetry}; see Section \ref{SSSLA}. 

The generalized notions are called with ``super" attached, in principle. 
For example, 
a $k$-\emph{super-functor} is a set-valued functor $\X$ defined on the category
$\mathtt{SAlg}_k$ of the super-commutative super-algebras (over $k$).
Important is that given such an $\X$, there is  
associated the $k$-functor $\X_{\mathsf{ev}}$  
which is
obtained from $\X$, restricting it to $\mathtt{Alg}_k$;
ordinary algebras are regarded as \emph{purely even} super-algebras,
or namely, super-algebras consisting of even elements, only.
We remark that the $k$-functor denoted here by $\X_{\mathsf{ev}}$  is alternatively denote by
$\X_{\mathsf{res}}$ in \cite{MZ1}, in which the symbol $\X_{\mathsf{ev}}$ is used to
denote the $k$-super-functor $R \mapsto \X(R_0)$.

The category of $k$-super-functors includes the following full subcategories:
\begin{equation*}
(\mathsf{faisceaux})\supset (\mathsf{faisceaux\ dur})\supset (\mathsf{super\text{-}schemes})
\supset (\mathsf{affine\ super\text{-}schemes})
\end{equation*}
The definition of affine super-schemes as well as that of affine super-groups will be
obvious. 
A \emph{faisceaux dur} (resp., \emph{faisceaux}) is a $k$-super-functor which preserves
finite direct products and every equalizer diagram
\[
R \to S \rightrightarrows S \otimes_R S
\]
that naturally arises from an \emph{fpqc} (resp., \emph{fppf}) morphism $R \to S$ in 
$\mathtt{SAlg}_k$. Here, fpqc (resp., fppf) represents ``faithfully flat"
(resp., ``faithfully flat and finitely presented"); see Section \ref{SSFF}.
A \emph{super-scheme} is a $k$-super-functor $\X$ which is \emph{local} (or roughly speaking,
$\X$ behaves like a sheaf with respect to Zariski coverings $\{ R \to R_{f_i}\}$, 
where $f_i \in R_0$ (finitely many) with $\sum_i f_iR_0=R_0$), and which is 
covered by some open sub-functors that are isomorphic to
affine super-schemes. By convention every super-scheme $\X$ is supposed to be \emph{non-trivial}
in the sense that $\X(R)\ne \emptyset$ for some $0\ne R\in \mathtt{SAlg}_k$. 

If a $k$-super-functor $\X$ is contained in one of the categories noted above, 
then the associated 
$k$-functor $\X_{\mathsf{ev}}$
is in the corresponding one defined in the ordinary (non-super) context. If $\G$ is an affine super-group, then $\G_{\mathsf{ev}}$
is an affine group.

There is an alternative definition (see \cite[Sect. 3.3]{CCF}, \cite[Chapter 4, Sect.~1]{Manin} or \cite[Sect.~4]{MZ1}, and 
rely on Leites \cite{L} for historical background)
of an (affine) super-scheme, 
which defines it as a super-ringed space satisfying some conditions that are naturally inferred
from the ordinary setup. 
We will use 
in part (in fact, essentially only in Section \ref{STHMII})
the word in this alternative, geometrical meaning. 
This is justified
since we have a natural equivalence (see \cite[Theorem 5.14]{MZ1}) between
the two categories of (affine) super-schemes thus defined in the two distinct ways. 
The equivalence of one direction 
assigns
to every super-scheme $\mathsf{Z}$ in the geometrical sense,
the $k$-super-functor 
$
R \, \mapsto \, \operatorname{Mor}(\operatorname{Spec} R,\, \mathsf{Z})
$
represented by it.

\subsection{The first main theorem}\label{SSMainI}
Suppose that an affine super-group $\G$ acts on a super-scheme $\X$. Here and in what follows actions are supposed to be on the right. 
We say that the action is \emph{free}, or $\G$ acts
\emph{freely} on $\X$, 
if the morphism of $k$-super-functors
\begin{equation}\label{EfALPHA}
\X \times \G \to \X \times \X,\quad (x,g) \mapsto (x,xg)
\end{equation}
is a monomorphism of super-schemes. 
Given a (free) $\G$-action on $\X$, there
is naturally induced a (free) action by the affine group 
$\G_{\mathsf{ev}}$
associated with $\G$ on the scheme $\X_{\mathsf{ev}}$ associated with $\X$. 

Suppose that an affine super-group $\G$ acts freely on a super-scheme $\X$.

\begin{definition}\label{DTors}
A morphism $\X \to \Y$ of super-schemes is called a \emph{super 
$\G$-torsor} (in the fpqc
topology), if it constitutes the co-equalizer diagram of faisceaux dur
\begin{equation}\label{ECoEQ}
\X \times \G \rightrightarrows \X \to \Y,
\end{equation}
where the paired arrows denote the action and the projection.
Alternatively, we say that $\X$ is a \emph{super $\G$-torsor over} 
$\Y$. 

A morphism
of super $\G$-torsors over $\Y$ is a $\G$-equivariant morphism of super-schemes over $\Y$;
such is necessarily an isomorphism, as will be seen in the paragraph following
Proposition \ref{PXG}. 
\end{definition}

\begin{rem}\label{RDG}
Demazure and Gabriel \cite[III, \S 4, 1.3 D\'{e}finition, p.361]{DG} define the notion of 
torsors in the category of faisceaux, where the relevant objects are
faisceaux and faisceaux en groupes. 
Our definition above, reduced to the
ordinary situation, becomes equivalent to their definition, when the
relevant 
objects are supposed to be schemes and affine algebraic groups; this fact
is easily seen from \cite[III, \S 4, 1.7 Corollaire, p.362]{DG} and Corollary \ref{CXG} below. 
\end{rem}

An affine super-group $\G$ is uniquely represented by a Hopf super-algebra.  
If the Hopf super-algebra is finitely generated, $\G$ is called an
\emph{affine algebraic super-group}. 

\begin{example}\label{EGoverH}
Let $\G$ be an affine algebraic super-group, and let $\mathsf{H}$ be a closed 
super-subgroup of $\G$. Then $\mathsf{H}$ acts freely on $\G$ by right-multiplication.
It is proved by \cite[Theorem 0.1]{MZ1}
 (and re-proved by \cite[Theorem 4.12]{MT}) that we have a Noetherian
(in fact, algebraic) super-scheme $\G/\mathsf{H}$ which fits into the co-equalizer diagram 
$\G \times \mathsf{H} \rightrightarrows \G \to \G/\mathsf{H}$ of faisceaux (dur); 
therefore, $\G/\mathsf{H}$ may be presented alternatively by 
$\G\tilde{/}\mathsf{H}$ (and by $\G\tilde{\tilde{/}}\mathsf{H}$), by using the notation 
which will be introduced in the following subsection.
One concludes that $\G\to \G/\mathsf{H}$ is a super $\mathsf{H}$-torsor. 
\end{example}

The first main theorem is the following; it shows, under some
restricted situation that involves affinity assumption, a remarkable
one-to-one correspondence between super-torsors and 
ordinary torsors.

\begin{theorem}\label{TBij}
Let $\G$ be a smooth affine algebraic super-group, and let $\Y$ be a Noetherian smooth
affine super-scheme. If $\X \to \Y$ is a super $\G$-torsor, then 
the induced morphism $\X_{\mathsf{ev}}\to \Y_{\mathsf{ev}}$ of schemes is 
a $\G_{\mathsf{ev}}$-torsor.
Moreover, 
the assignment $\X \mapsto \X_{\mathsf{ev}}$ gives rise to a bijection
from
\begin{itemize}
\item the set of all isomorphism classes of super $\G$-torsors over $\Y$
\end{itemize}
onto
\begin{itemize}
\item the set of all isomorphism classes of $\G_{\mathsf{ev}}$-torsors over 
$\Y_{\mathsf{ev}}$.
\end{itemize}
\end{theorem}

\begin{rem}\label{RANS}
For those super $\G$-torsors $\X\to \Y$ and 
$\G_{\mathsf{ev}}$-torsors 
$\mathsf{Z}\to \Y_{\mathsf{ev}}$ which
constitute the two sets above, 
$\X$ and $\mathsf{Z}$ are necessarily 
affine, Noetherian and smooth. 
The first two properties follow from Proposition \ref{PXG} (2) below.
The remaining smoothness will be proved by  Proposition \ref{PHS}.
\end{rem}

Since in Theorem \ref{TBij}, 
everything thus turns out to be affine, the theorem is re-formulated in terms of
(Hopf) super-algebras which represent the relevant functors. 
A super $\G$-torsor 
$\X \to \Y$ with $\X$, $\Y$ affine is translated into an $H$-\emph{Galois extension} $A/B$
(see Definition \ref{DGAL}),
where $H$ is the Hopf super-algebra representing $\G$, and $A$ and $B$ are the super-algebras
representing $\X$ and $\Y$, respectively.
Theorem \ref{TBij} will be obtained as an immediate consequence
of Theorem \ref{THA}, which explicitly describes such a Hopf-Galois extension
in the super situation, in terms of the naturally associated 
extension in the ordinary situation; see also Corollary \ref{CHA}. 

\begin{rem}\label{RNoHope}
(1)\
It will be seen in Remark \ref{RNoncat} that 
the bijection shown in Theorem \ref{TBij} does not
result from a category-equivalence from 
the category of super $\G$-torsors over $\Y$ to
the category of $\G_{\mathsf{ev}}$-torsors over $\Y_{\mathsf{ev}}$; these categories are in fact
groupoids, as is seen from the paragraph following Definition \ref{DTors}.

(2)\
A super-scheme $\X$, regarded as a super-ringed space, is said to be \emph{split} (see \cite[Section 2.6]{MT}, for example), if 
the structure sheaf $\cO_{\X}$ is isomorphic to the exterior algebra 
$\wedge_{\cO_{\X_{\mathsf{ev}}}}(\mathcal{M})$ on some locally free 
$\cO_{\X_{\mathsf{ev}}}$-module sheaf $\mathcal{M}$, where $\cO_{\X_{\mathsf{ev}}}$ denotes
the structure sheaf of the scheme associated with $\X$; see the first paragraph of Section \ref{STHMII}.
Every affine algebraic super-group $\G$ is split (see \eqref{ETPD}--\eqref{ETPDop} below), 
while $\G/\mathsf{H}$ such as in Example \ref{EGoverH} 
can be non-split; see \cite[Remark 4.20 (2)]{MT}. 
In view of the cited example, 
it would be hopeless to generalize Theorem \ref{TBij}, 
removing the affinity assumption. 
\end{rem}

\subsection{The second main theorem}\label{SSMainII}
Suppose that an affine super-group $\G$ acts freely on a super-scheme $\X$. We 
have the $k$-super-functor which assigns to each $R \in \mathtt{SAlg}_k$, 
the set $\X(R)/\G(R)$ of $\G(R)$-orbits; it has 
the property 
\begin{equation}\label{ENiceP}
\X(R)/\G(R)\subset \X(S)/\G(S)\ \, \text{whenever}\ \, R \subset S,
\end{equation}
since
the action is supposed to be free. We can construct in a simple 
manner (see \cite[Remark 3.8]{MZ1}, for example)
the faisceau dur $\X\tilde{\tilde{/}}\G$ which includes the $k$-super-functor above
as a sub-functor, and which 
is \emph{associated with} 
that $k$-super-functor in the following sense:
every morphism from the $k$-super-functor
to any faisceau dur extends uniquely to a morphism from $\X\tilde{\tilde{/}}\G$. Clearly,
this $\X\tilde{\tilde{/}}\G$ is characterized by the co-equalizer diagram of faisceaux dur
$\X \times \G \rightrightarrows \X \to \X\tilde{\tilde{/}}\G$ which is analogous to \eqref{ECoEQ}. An important consequence is the following.

\begin{prop}\label{PXG}
Let $\G$ be an affine super-group. 
\begin{itemize}
\item[(1)]
A super $\G$-torsor refers precisely 
to the natural morphism $\X \to \X\tilde{\tilde{/}} \G$, where $\G$ acts freely on a super-scheme $\X$ and the associated faisceau dur $\X \tilde{\tilde{/}}\G$ is a super-scheme. 
\textup{(}Here, strictly speaking, one identifies a super $\G$-torsor $\X \to \Y$ with
$\X \to \X \tilde{\tilde{/}}\G$ such as above, through a unique isomorphism 
$\Y \overset{\simeq}{\longrightarrow}
\X \tilde{\tilde{/}} \G$ 
compatible with the morphisms from $\X$.\textup{)}
\item[(2)]
Suppose that $\X \to \Y$ is a super $\G$-torsor. Then it is an affine and faithfully flat
morphism of super-schemes. 
Moreover,
\begin{equation}\label{EXGXYX}
\X \times \G \to \X \times_{\Y}\X,\quad (x,g)\mapsto (x, xg)
\end{equation}
is an isomorphism of super-schemes. If $\G$ is an affine algebraic super-group, then
$\X \to \Y$ is a morphism of finite presentation. 
\end{itemize}
\end{prop}

This is proved essentially by the same argument as proving \cite[Part I, 5.7, (1)]{J}, 
which, in fact, uses 
the property \eqref{ENiceP}. 
The faithful flatness shown in Part 1 above implies that
a morphism of super $\G$-torsors over some $\Y$ is necessarily an isomorphism. 

Note that the last assertion of Part 2 above follows from
the algebraicity assumption for $\G$, by using \eqref{EXGXYX}, base-change
and faithfully-flat descent. 
As this assertion suggests, 
it would be more natural to consider the faisceau $\X\tilde{/}\G$
which is associated, now in 
an obvious sense, 
with the $k$-super-functor $R \mapsto \X(R)/\G(R)$, rather than the faisceau dur
$\X\tilde{\tilde{/}}\G$ above, when we assume that $\G$ is algebraic.
In general (without the algebraicity assumption), one sees from the constructions that
\[
(\X \tilde{/}\G)(R) \subset (\X\tilde{\tilde{/}}\G)(R)
\]
for every $R \in \mathtt{SAlg}_k$. 

From the above-mentioned assertion one sees the following.

\begin{corollary}\label{CXG}
Suppose that an affine algebraic super-group $\G$ acts freely on a super-scheme $\X$.
Then the faisceau $\X \tilde{/}\G$
is a super-scheme if and only if the faisceau dur
$\X \tilde{\tilde{/}} \G$ is.
If these equivalent conditions are satisfied, 
we have 
\[
\X\tilde{/}\G= \X\tilde{\tilde{/}}\G,
\]
and $\X \to \X\tilde{/}\G\, (= \X\tilde{\tilde{/}}\G)$ is a 
super $\G$-torsor. 
\end{corollary}

The second main theorem is the following; it shows, under some assumptions
that include (i)~$\G$ is algebraic and (ii)~the $\G$-action is free, that  
$\X\tilde{/}\G$ is a super-scheme provided $\X_{\mathsf{ev}}\tilde{/}\G_{\mathsf{ev}}$ 
is a scheme. Due to (i), we thus work with faisceaux, 
not with faisceaux dur.

\begin{theorem}\label{THMII}
Suppose that an affine algebraic super-group $\G$ acts freely
on a locally Noetherian super-scheme $\X$. Assume the following:
\begin{itemize}
\item[(a)] $\G$ and $\X$ are both smooth; 
\item[(b)] The faisceau $\X_{\mathsf{ev}}\tilde{/}\G_{\mathsf{ev}}$ associated with the
induced \textup{(}necessarily, free\textup{)} action by the affine algebraic group 
$\G_{\mathsf{ev}}$ on the locally Noetherian scheme $\X_{\mathsf{ev}}$ 
is a scheme; it then necessarily holds that $\X_{\mathsf{ev}}\tilde{/}\G_{\mathsf{ev}}$ is locally Noetherian and smooth, and 
$\X_{\mathsf{ev}}\to \X_{\mathsf{ev}}\tilde{/}\G_{\mathsf{ev}}$ is 
a $\G_{\mathsf{ev}}$-torsor. 
\end{itemize}
Then the faisceau $\X \tilde{/} \G$ is a locally Noetherian, smooth super-scheme such that
\begin{equation}\label{EXGXG}
\X_{\mathsf{ev}}\tilde{/}\G_{\mathsf{ev}}=(\X \tilde{/} \G)_{\mathsf{ev}}. 
\end{equation}
In particular, $\X\to \X\tilde{/}\G$ is a super $\G$-torsor.
\end{theorem}

By \eqref{EXGXG} we mean that the canonical morphism 
$\X_{\mathsf{ev}}\to (\X \tilde{/} \G)_{\mathsf{ev}}$
induces an isomorphism such as noted.
By convention we thus present canonical isomorphisms by $=$, 
the equal sign. 

To prove the theorem above, essential is the case when $\X$ is affine, in which case the theorem is re-formulated,
again in terms of Hopf-Galois extensions,
by Theorem \ref{THB}.

\begin{rem}
One might think that the result \cite[Theorem 4.12]{MT} referred to in Example \ref{EGoverH} follows from
the above Theorem \ref{THMII} as the special case where $\mathsf{H}$ and $\G$ of the former are taken as 
$\G$ and $\X$, respectively, of the latter, 
since the faisceau $\G_{\mathsf{ev}}\tilde{/}\mathsf{H}_{\mathsf{ev}}$ is known 
to be a scheme (see \cite[Part I, Sect.~5.6, (8)]{J}), and thus the assumption (b)
above is satisfied. But this is true only when the assumption (a) (for $\mathsf{H}$ and $\G$) as well is satisfied,
which is always the case if $\operatorname{char}k=0$; see Remark \ref{RSmooth} (1). 
\end{rem}
\subsection{Organization of the paper and additional remarks}\label{SSOrg}
As was mentioned above, the two main theorems, Theorems \ref{TBij} and \ref{THMII}, 
largely depend on Theorems \ref{THA} and
\ref{THB}, which are formulated in Section \ref{SHA} in terms of Hopf-Galois extensions,
the algebraic counterpart of super-torsors. Those Hopf-algebraic theorems 
are proved by spending the whole of Section \ref{STHA} and of Section \ref{STHB}, respectively. 
Theorem \ref{THMII} is proved in the final Section \ref{STHMII}, while
Theorem \ref{TBij} is an immediate consequence of Theorem \ref{THA},
as was already mentioned.
Theorem \ref{THA}
is accompanied with an analogous result, Proposition \ref{PHA}, which 
is used to prove Theorem \ref{THB};
Section \ref{SPHA} is devoted to proving the proposition. 
Section \ref{SSP} is devoted to preliminaries for 
the following Sections \ref{SHA}--\ref{STHB}. 

Sections \ref{SSP}--\ref{STHB} will be found to have a purely Hopf-algebraic flavor. 
We there use the Hopf-algebraic techniques developed so far to investigate
algebraic super-groups, as well as several results already obtained by the techniques;
see \cite{M}, \cite{M1}, \cite{MO}--\cite{MZ2}. 
In addition, there are applied some important results proved by 
Schauenburg and Schneider \cite{SS} for 
Hopf-Galois extensions in the ordinary but non-commutative setting; see Sections 
\ref{SSHG} and \ref{STHA1}. A notable result of the application is
the equivariant smoothness, Proposition \ref{PHS}, proved in the latter Section \ref{STHA1}. 

\begin{rem}\label{RContinue} 
The two main theorems of ours both play substantial roles in the preprint
\cite{Masuoka} (continuing to \cite{MasuokaII}) 
by the first-named author; it generalizes Picard-Vessiot theory, which is Galois theory of 
linear differential equations, to the super context. 
\end{rem}

\section{Preliminaries}\label{SSP}

\subsection{Super-linear algebra}\label{SSSLA}
We let $\ot$ denote the tensor product $\ot_k$ over our base field $k$;
recall $\operatorname{char} k\ne 2$ by assumption.
A vector space $V= V_0\oplus V_1$ graded by the group $\mathbb{Z}_2=\{ 0,1\}$
is alternatively called a \emph{super-vector space}; 
it is said to be \emph{purely even} (resp., \emph{purely odd})
if $V=V_0$ (resp., if $V=V_1$). The super-vector spaces form a
symmetric tensor category,
\[
\M_k=(\M_k,\ot,k),
\]
with respect to the so-called super-symmetry
\[
c=c_{V,W}:V\ot W\overset{\simeq}{\longrightarrow} W\ot V, \quad c(v\ot w)=(-1)^{|v||w|}
w\ot v,
\]
where $V,W \in \M_k$, and $v\, (\in V)$ and $w\, (\in W)$ are homogeneous
elements of degrees $|v|$, $|w|$. 
A \emph{super-algebra} is by definition an algebra in $\M_k$. Similarly the notions of
\emph{super-coalgebra}, \emph{Hopf super-algebra} and \emph{Lie super-algebra}
are defined in an obvious manner.

\subsection{Super-(co)algebras}\label{SSSCA}
Given a super-algebra $R$, we let
\[ \M_R \quad \textup{(resp.,}\  {}_R\M \textup{)} \]
denote the category of right (resp., left) $R$-modules in $\M_k$; an object in
the category is called a right (resp., left) $R$-super-module. As is shown in
\cite[Lemma 5.1 (2)]{M}, 
such a super-module is projective in the category if and only if
it is projective, regarded as an ordinary $R$-module. 

Given a super-coalgebra $C$, we let
\[ \M^C \quad \textup{(resp.,}\  {}^C\M \textup{)} \]
denote the category of right (resp., left) $C$-comodules in $\M_k$; an object in
the category is called a right (resp., left) $C$-\emph{super-comodule}. 
It follows by the dual argument of proving
\cite[Lemma 5.1 (2)]{M} referred to above that
such a super-comodule is injective in the category if and only if
it is injective, regarded as an ordinary $C$-comodule. Suppose that 
$V=(V,\rho)$ is an object
in $\M^C$ and $W=(W,\lambda)$ is in ${}^C\M$. The \emph{co-tensor product} $V\square_CW$
is defined by the equalizer diagram
\begin{equation}\label{ECoTENS}
V\square_CW \to V\otimes W \rightrightarrows V\otimes C\otimes W
\end{equation}
in $\M_k$, 
where the paired arrows indicate
$\rho \ot \mathrm{id}_W$ and $\mathrm{id}_V\otimes \lambda$. This gives rise to 
the functor $V\square_C: {}^C\M\to \M_k$, which is seen to be left exact. 
It is known that  
the functor $V\square_C$
is exact if and only if $V$ is injective. An analogous result for
the analogous functor $\square_CW$ holds. 

For $V =(V,\rho)$ as above, we define
\begin{equation}\label{ErCinv}
V^{\mathrm{co}\hspace{0.2mm}C}:=\{\, v \in V\mid \, \rho(v)=v \ot 1\, \}.
\end{equation}
This is a super-subvector space of $V$, whose elements are called 
$C$-\emph{co-invariants}
in $V$. Analogously, for $(W,\lambda)$ as above, we define
\begin{equation}\label{ElCinv}
{}^{\mathrm{co}\hspace{0.2mm} C}W:=\{\, w\in W \mid \lambda(w)=1\ot w \, \}.
\end{equation}

In what follows we assume that all super-algebras $R$ are \emph{super-commutative},
or namely, $ba=(-1)^{|a||b|}ab$ for all homogeneous elements $a, b$; the assumption 
is equivalent to saying that $R$ includes the even component $R_0$ as a central
subalgebra, and we have $ba=-ab$ for all $a , b \in R_1$. 
The two categories 
of super-modules are then naturally identified so that
\[
{}_R\M=\M_R. 
\]

\subsection{Hopf super-algebras}\label{SSHSA}
Accordingly, all Hopf super-algebras $H$ are assumed to be super-commutative,
unless otherwise stated. The category
\[ \M^H =(\M^H, \, \ot,\ k) \]
of right $H$-super-comodules then forms a symmetric tensor category with respect to
the super-symmetry. 
We will mostly deal 
with these right $H$-super-comodules rather than the analogous left ones.

In general, given a super-algebra $R$, we let
\begin{equation}\label{EBAR}
I_R:=(R_1),\quad \overline{R}:=R/I_R
\end{equation}
denote the super-ideal of $R$ generated by the odd component $R_1$, and the
quotient super-algebra of $R$ by $I_R$, respectively; the latter is in fact the largest purely-even quotient super-algebra of $R$. 

Let $H$ be a Hopf super-algebra. Then $\overline{H}$ is seen to be an ordinary Hopf 
algebra.
Let $H^+$ be the augmentation super-ideal of $H$, or namely, the kernel
$\mathrm{Ker}(\varepsilon :H \to k)$ of the co-unit of $H$. Let
\begin{equation}\label{EW}
W:=(H^+/(H^+)^2)_1 
\end{equation}
denote the odd component of the super-vector space $H^+/(H^+)^2$. 
The exterior algebra $\wedge(W)=\bigoplus_{n\ge 0}\wedge^n(W)$ on $W$ is regarded 
as a super-algebra whose even (resp., odd) component is the direct sum of those $\wedge^n(W)$ with $n$ 
even (resp., odd).
Since
$H$ is naturally an algebra in $\M^{H}$, it turns to be an
algebra in $\M^{\overline{H}}$ along the quotient map $H \to \overline{H}$. 
The tensor-product decomposition theorem \cite[Theorem 4.5]{M} 
tells us that that there is a co-unit-preserving isomorphism
\begin{equation}\label{ETPD}
H \simeq \wedge(W)\ot \overline{H}
\end{equation}
of algebras in $\M^{\overline{H}}$, where the co-unit of the right-hand side
is the co-unit of $\overline{H}$ tensored with 
the natural projection $\wedge(W) \to \wedge^0(W)=k$. 
As an exception from what we said above about the side, we will mainly use the 
opposite-sided version 
\begin{equation}\label{ETPDop}
H \simeq \overline{H}\ot \wedge(W)
\end{equation}
of \eqref{ETPD}.

\begin{rem}
(1)\ 
As a typical example of affine algebraic super-groups, recall $\mathsf{GL}_{m|n}$; for a super-vector 
space $V$ with $m=\dim V_0$, $n=\dim V_1$, the super-group assigns
to every super-algebra $R$, the group of all automorphisms of the $R$-super-module $R\otimes V$.
The Hopf super-algebra representing $\mathsf{GL}_{m|n}$ is presented, for example, in Example 4.1 of 
\cite{M2}. 
Example 4.2 of \cite{M2} explicitly gives a relevant isomorphism such as \eqref{ETPDop}; see Footnote
5 added to the arXiv version, in particular.

(2)\
We emphasize that the isomorphisms \eqref{ETPD} and \eqref{ETPDop} preserve
the $\overline{H}$-co-actions, in particular, and they hold without assuming that $H$ (or $\overline{H}$) 
is finitely generated or smooth. They were first proved by \cite[Proposition 2.4]{MO} under that
assumption that $\operatorname{char}k=0$ (whence $\overline{H}$ is smooth), and $H$ is finitely
generated; we remark that the cited proposition refers to a result by H.~Boseck as a weaker version,
but his proof was wrong.
\end{rem}

\subsection{Smooth super-algebras}\label{SSSSA}
A super-algebra $R$ is said to be \emph{Noetherian}, if the super-ideals of $R$ satisfy
the ACC, or equivalently, if the even component $R_0$ is a Noetherian algebra, and 
the odd component $R_1$, regarded as an $R_0$-module,
is finitely generated; see \cite[Section 2.4]{MT} 
for other equivalent conditions. 

A super-algebra $R$ is said to be \emph{smooth} (over $k$), 
if given a super-algebra surjection $A \to B$ with nilpotent kernel, every super-algebra map 
$R\to B$ factors through $A$ by some super-algebra map $R \to A$. 
It is known that if $\operatorname{char}k=0$, every 
Hopf super-algebra is smooth; see Remark \ref{RSmooth} (1) below.

The following is part of \cite[Theorem A.2]{MZ2}.

\begin{prop}\label{PSMOO}
For a Noetherian super-algebra $R$, the following conditions \textup{(a)} and \textup{(b)}
are equivalent:
\begin{itemize}
\item[(a)] $R$ is smooth as a super-algebra;
\item[(b)] 
\begin{itemize}
\item[(i)] $\overline{R}$ is smooth as an algebra,
\item[(ii)] $I_R/I_R^2$, which is seen to be a purely odd $\overline{R}$-super-module, is finitely generated projective as an $\overline{R}$-module, and 
\item[(iii)] there is an isomorphism $\wedge_{\overline{R}}(I_R/I_R^2)\simeq R$ of super-algebras.
\end{itemize}
\end{itemize}
\end{prop}

In (iii) above, 
$\wedge_{\overline{R}}(I_R/I_R^2)$ denotes the exterior $\overline{R}$-algebra on $I_R/I_R^2$, which 
is regarded as a super-algebra 
just as the $\wedge(W)$ before. One can replace $\overline{R}$ and $I_R/I_R^2$, more generally,
with a purely even super-algebra and a purely odd super-module over it, respectively. 

\begin{rem}\label{RSMOO}
(1)\
For a super-algebra $R$ in general, one sees as in \cite[Page 360, line --9]{MZ2} that 
\[
I_R/I_R^2= R_1/R_1^3,
\]
which is a purely odd $\overline{R}$-super-module. This is finitely generated,
if $R$ is Noetherian; recall that $R_1$ is then finitely generated as an $R_0$-module.

(2)\ Suppose that a Noetherian super-algebra $R$ satisfies the equivalent conditions
in the preceding proposition. By (b)(i), there is a section of the natural projection
$R_0 \to \overline{R}\, (=R_0/R_1^2)$, through which we regard $R$ as a super-algebra over 
$\overline{R}$. The natural projection
$I_R \to I_R/I_R^2$ is seen to be a surjection in $\M_{\overline{R}}$, which has a
section
by (b)(ii). The proof of \cite[Theorem A.2]{MZ2} shows that
the $\overline{R}$-super-algebra map
\begin{equation}\label{EWIR}
\wedge_{\overline{R}}(I_R/I_R^2)\to R
\end{equation}
which uniquely extends the last section
is an isomorphism such as in (b)(iii).
\end{rem}

\subsection{Faithful flatness and finite presentation}\label{SSFF}
A map $R \to S$ of super-algebras is said to be (\emph{faithfully}) \emph{flat}
if the functor $S\ot_R : {}_R\M \to \M_k$ is (faithfully) exact. 
This is equivalent to saying that $S$ is (faithfully) flat, regarded as an ordinary 
right $R$-module; see \cite[Lemma 5.1 (1)]{M1}. Each of the equivalent conditions is equivalent to the
one with the side switched. 

Given a map $R \to S$ of super-algebras, 
$S$ is said to be \emph{finitely presented} over $R$, if $S$ is presented 
so as $R[x_1,\cdots,x_m;y_1,\cdots,y_n]/I$, where $x_1,\dots,x_m$ are finitely many even variables, $y_1,\dots,y_n$
are finitely many odd variables, and $I$ is a super-ideal which is generated by finitely many homogeneous elements. 

Suppose that $\X$ is an affine super-scheme. 
We let $\cO(\X)$ denote the super-algebra representing 
$\X$; it is non-zero since $\X$ is 
non-trivial by convention; see Section \ref{SSMiP}. 
Note that the algebra $\cO(\X_{\mathsf{ev}})$ 
which represents the associated affine scheme $\X_{\mathsf{ev}}$
equals $\overline{\cO(\X)} \, (=\cO(\X)/(\cO(\X)_1))$.
We say that $\X$ \emph{has a property} (P) (e.g., Noetherian,
smooth), if $\cO(\X)$ has the property (P).

\section{Shifting to the Hopf-algebra world}\label{SHA}

\subsection{Hopf-Galois extensions}\label{SSHG}

Let $\G$ be an affine super-group, and set 
$H:=\cO(\G)$, the Hopf super-algebra 
representing $\G$. A left $\G$-super-module may be understood to be
a right $H$-super-comodule. 

Let $\X$ is an affine super-scheme, and set $A:=\cO(\X)\, (\ne 0)$.
Suppose that $\G$ acts on $\X$ from the right, and the action 
$\X\times \G \to \X$ is represented by 
\begin{equation}\label{ERHO}
\rho : A \to A \ot H. 
\end{equation}
Thus, $A=(A, \rho)$ is an algebra in $\mathtt{SMod}^H$. Let 
\begin{equation}\label{ECOI}
B :=A^{\mathrm{co}\hspace{0.2mm} H}\, (=\{ a \in A\mid \rho(a)=a\otimes 1\});
\end{equation}
see \eqref{ErCinv}. This is a super-subalgebra of $A$, and is often
alternatively presented as $B=A^G$, the $G$-invariants in $A$. 

Recall that the $\G$-action on $\X$ is said to be
\emph{free}, if the morphism \eqref{EfALPHA}
of $k$-super-functors 
is a monomorphism of super-schemes, or equivalently,
if the natural maps
\[ \X(R) \times \G(R) \to \X(R) \times \X(R),\ (x,g)\mapsto (x, xg) \]
are injective for all $R \in \mathtt{SAlg}_k$. This is the case if
the super-algebra map
\begin{equation}\label{EALP}
\alpha : A \otimes A \to A \otimes H,\quad \alpha(a \ot b)=a\, \rho(b), 
\end{equation} 
which represents \eqref{EfALPHA},
is surjective. Notice that this super-algebra map 
induces
\begin{equation}\label{EBET}
\beta : A \otimes_B A \to A \otimes H,\quad \beta(a \ot_B b)=a\, \rho(b). 
\end{equation}

\begin{definition}\label{DGAL}
We say that $A/B$ is (an) $H$-\emph{Galois} (\emph{extension}), 
or $A$ is an $H$-\emph{Galois extension
over} $B$, if $B \hookrightarrow A$ is faithfully flat, and the map $\beta$ above is bijective. 
\end{definition}

\begin{rem}\label{RHG}
The definition above is generalized to the situation where $H$ and $A$ are not necessarily
super-commutative. But it seems reasonable to assume then that the antipode of $H$ is bijective,
in which case, the condition that $A$ is faithfully flat over $B$ on the left is equivalent to
the condition that it is so on the right, as is seen from 
\cite[Theorem~4.10]{SS}; see also \cite[Lemma 10.1]{MZ1} and Remark \ref{RGUE} below. 
\end{rem}

Recall that $A$ is an algebra in $\M^H$. We let
\begin{equation}\label{ESMod}
\mathtt{SMod}_A^H= (\mathtt{SMod}^H)_A
\end{equation}
denote the category of right $A$-modules in $\mathtt{SMod}^H$. Its objects may be called
$(H,A)$-\emph{Hopf super-modules} in view of \cite[p.244]{Doi}. 
One sees $A \in \mathtt{SMod}_A^H$. More generally, $N\ot_BA\in \mathtt{SMod}_A^H$, if   
$N\in \mathtt{SMod}_B$. In fact we have the functor
\begin{equation}\label{EPSI}
\Psi : \mathtt{SMod}_B\to \mathtt{SMod}_A^H,\quad \Psi(N)=N\ot_B A.
\end{equation}
This is left adjoint to the functor $\Phi: \mathtt{SMod}_A^H\to \mathtt{SMod}_B$
which assigns to each $M$ in $\mathtt{SMod}_A^H$,
the $H$-co-invariants $M^{\mathrm{co}\hspace{0.2mm} H}$ in $M$; see \eqref{ErCinv}. 

We reproduce the following from \cite{MZ1}. See also \cite[Sect. 4]{Zubkov}, \cite[Remark 2.8]{MSS}. 

\begin{theorem}[\text{\cite[Theorem 7.1]{MZ1}}]\label{TOBER}
Retain the situation as above.
\begin{itemize}
\item[(1)] The following are equivalent: 
\begin{itemize}
\item[(i)]
The $\G$-action on $\X$ is free, and the faisceau dur $\X \tilde{\tilde{/}} \G$ is 
an affine super-scheme;
\item[(ii)]
The map $\alpha$ is surjective, and 
$A$ is injective in $\mathtt{SMod}^H$;
\item[(iii)]
$A/B$ is an $H$-Galois extension;
\item[(iv)]
The functor $\Psi$ in \eqref{EPSI} is an equivalence.   
\end{itemize}
If these equivalent conditions are satisfied, then the affine super-scheme
$\X \tilde{\tilde{/}} \G$ 
is represented by $B$, or in notation, $\cO(\X \tilde{\tilde{/}} \G)=B$.
\item[(2)]
Suppose that $\G$ is algebraic and $\X$ is Noetherian \textup{(}or in other words, $H$ is finitely generated and $A$ is Noetherian\textup{)}. 
If the equivalent conditions above are satisfied, 
then $\X \tilde{\tilde{/}} \G$ \textup{(}or $B$\textup{)} is Noetherian, 
and it coincides with the faisceau $\X \tilde{/} \G$.
\end{itemize}
\end{theorem}

To prove part of this, the same argument as we used above to see Proposition \ref{PXG}
had been used. Notice from Corollary \ref{CXG} 
that the conclusion $\X \tilde{\tilde{/}} \G=\X \tilde{/} \G$
in Part 2 holds, only assuming that $\G$ is algebraic. Part 1 of the theorem above,
combined with Part 2 of Proposition \ref{PXG}, shows the following. 

\begin{corollary}\label{CXGAB}
Let $\G$ and $H$ be as above. Let $\Y$ be an affine super-scheme, and set $B:=\cO(\Y)$. 
Every super $\G$-torsor
$\X$ over $\Y$ is necessarily an affine super-scheme, and $\cO(\X)$ naturally turns
into an $H$-Galois extension over $B$. Moreover, every $H$-Galois
extension over $B$ arises uniquely from a super $\G$-torsor over $\Y$ in this way. 
\end{corollary}

The following is well known at least in the ordinary situation, and is easily proved
by using faithfully flatness of Hopf-Galois extension; see also the paragraph
following Proposition \ref{PXG}. 

\begin{lemma}\label{LIsom}
Let $H$ be a Hopf super-algebra, and let $B$ be a non-zero super-algebra.
Given $H$-Galois extensions $A/B$ and $A'/B$, every $B$-algebra morphism $A \to A'$ in $\M^H$
is necessarily an isomorphism. 
\end{lemma}

We say that $A/B$ and $A'/B$ are \emph{isomorphic} to each other, if there exists 
a morphism, necessarily an isomorphism, such as above. 

The next Proposition \ref{PSS} shows an interesting property of $H$-Galois extensions, which will be used
to prove Lemmas \ref{LGF} and \ref{LSP}. In fact, it directly generalizes an important result  
by Schauenburg and Schneider \cite{SS}, which was proved in the non-super situation, to our super situation.

To prove the theorem above reproduced from \cite[Section 10]{MZ1}, there was used 
the so-called bosonization technique, in order to reduce
(main part of) the theorem 
to a known result which was
proved in the non-super, but non-commutative situation. 
The same technique is used, as will be seen soon, to prove the next proposition.
We continue to assume the super-commutativity.
But the proof works for proving the result without the
assumption; see Remark \ref{RHG}. 

Let $A/B$ be an $H$-Galois extension. Since $B$ is an algebra in $\M^H$ on which
$H$ co-acts trivially, the category
\[ {}_B\M^H\, (={}_B(\M^H)) \]
of left $B$-modules in $\M^H$ consists 
precisely
of those left $B$-super-modules $M$ which
are at the same time right $H$-super-comodules such that the $B$-action 
$B\ot M \to M$ is $H$-co-linear.

\begin{prop}\label{PSS}
Let $A/B$ be as above. Then the product map
\begin{equation}\label{EPR}
B \ot A \to A,\quad b\ot a \mapsto ba
\end{equation}
is clearly 
a surjective morphism in ${}_B\M^H$, and in fact, it splits as such a 
morphism. 
Hence, $A$ is projective in ${}_B\M$, in particular. 
\end{prop}
\pf
The bosonization of $H$ is the (ordinary) Hopf algebra
\[
\widehat{H}=H \lboson \mathbb{Z}_2
\]
of smash product (resp., smash co-product) by $\mathbb{Z}_2$
as an algebra (resp., as a coalgebra). 
This Hopf algebra
may not be commutative, but has a bijective antipode; see \cite[Lemma 10.1]{MZ1}.
In this proof we present $\mathbb{Z}_2$ as a multiplicative group
with generator $\sigma$. 
One sees from \cite[Proposition 10.3]{MZ1}
that the semi-direct product of $A$ by $\mathbb{Z}_2$
\begin{equation}\label{EAhat}
\widehat{A}=A \rtimes \mathbb{Z}_2
\end{equation}
naturally turns into an $\widehat{H}$-comodule algebra, 
and moreover, into an $\widehat{H}$-Galois
extension over $B=B\times \{ 1 \}$ (in the non-super situation), 
which is faithfully flat over $B$ on both sides. 
It follows by \cite[Theorem 5.6]{SS} that the product map
$B \ot \widehat{A}\to \widehat{A}$ has 
a left $B$-linear and right $\widehat{H}$-co-linear section,
\[
s : \widehat{A} \to B \ot \widehat{A}.
\]
We suppose $\widehat{A}=A\ot k\mathbb{Z}_2$, and define maps 
$r : B  \ot A\ot k\mathbb{Z}_2\to B\ot A$ and
$t : A \to A \ot k\mathbb{Z}_2$
by
\[
r(b \ot a \ot \sigma^i)= b\ot a,\ \, i=0,1;\quad t(a)= a\ot e, 
\]
where $e=\frac{1}{2}(1+\sigma)$, a primitive idempotent in $k\mathbb{Z}_2$. 
Then the composite $r\circ s\circ t: A\to B\ot A$
is seen to be a left $B$-linear and right $H$-co-linear section 
of the product map \eqref{EPR}, 
which may not preserve the parity. The section as well as the
product map is naturally regarded as a map of left $B \otimes H^*$-modules,
where $H^*$ is the dual algebra of $H$. Since
$B \otimes H^*\subset (B \otimes H^*)\rtimes \mathbb{Z}_2$ is a separable extension of rings by the assumption $\operatorname{char}k\ne 2$, 
one sees that the product map has a
$(B \otimes H^*)\rtimes \mathbb{Z}_2$-linear section, 
which is indeed a desired section; cf. the proof of \cite[Lemma 5.1 (2)]{M}. 
\epf

\subsection{Inflation of Hopf-Galois extensions}\label{SSINF}
G\"{u}nther \cite{G} proves in the non-super situation 
an interesting result on the inflation of Hopf-Galois
extensions, which is directly generalized to the super situation; 
see Remark \ref{RGUE} below. 
We are going to present the result in a restricted form as will be needed in the sequel. 
 
Continue to suppose that $\G$ is an affine super-group represented by $H$. 
Recall that the affine group $\G_{\mathsf{ev}}$ is represented by the quotient
Hopf algebra $\overline{H}=H/(H_1)$ of $H$. Let 
\[
\Delta : H \to H \ot H,\quad \Delta(h)=h_{(1)}\otimes h_{(2)}
\]
denote the co-product on $H$. Then $H$ turns into an algebra in 
${}^{\overline{H}}\mathtt{SMod}$ 
with respect to
$\overline{\Delta}:H \to \overline{H}\ot H$,\ 
$\overline{\Delta}(h)= \overline{h}_{(1)}\ot h_{(2)}$, 
where
$H \to \overline{H},\ h \mapsto \overline{h}$ denotes the quotient map. Let
\begin{equation}\label{ER}
R= {}^{\mathrm{co}\hspace{0.2mm} \overline{H}}H\, 
(=\{\, h \in H\, \mid \, \overline{\Delta}(h)=\overline{1} \otimes h \, \})
\end{equation}
denote the super-subalgebra of 
\emph{left} $\overline{H}$-co-invariants  (see \eqref{ElCinv}), which is seen to
satisfy $\Delta(R)\subset R\ot H$, whence it is an algebra in $\mathtt{SMod}^H$. 
From \eqref{ETPDop} we see
\begin{equation}\label{eqRwedgeW}
R \simeq \wedge(W).
\end{equation}

Suppose that $D=(D,\rho)$ is an algebra in $\mathtt{SMod}^{\overline{H}}$, 
which may not be purely even. The co-tensor product (see \eqref{ECoTENS})
\[
D \square_{\overline{H}} H=\Big\{ \, \sum_i d_i\ot h_i \in D\ot H\, \Big|
\, \sum_i \rho(d_i)\ot h_i=\sum_i d_i\ot \overline{\Delta}(h_i)\, \Big\}
\]
is naturally an algebra in $\mathtt{SMod}^H$ such that
\begin{equation}\label{EDcoH}
D^{\mathrm{co}\hspace{0.2mm} \overline{H}}
=(D \square_{\overline{H}} H)^{\mathrm{co}\hspace{0.2mm} H}.
\end{equation}
This canonical isomorphism \eqref{EDcoH} is a restriction of $D \to D \otimes H,\ d \mapsto d\ot 1$, and 
it has, indeed, an inverse, $\sum_i d_i\ot h_i\mapsto \sum_i d_i\varepsilon(h_i)$, given by the co-unit
$\varepsilon$ of $H$.

\begin{lemma}\label{LINF}
Set $B:=D^{\mathrm{co}\hspace{0.2mm} \overline{H}}$.
If $D/B$ is $\overline{H}$-Galois,
then $(D \square_{\overline{H}} H)/B$ is $H$-Galois. 
\end{lemma}

\pf
One sees from \eqref{EDcoH} that $B=(D \square_{\overline{H}} H)^{\mathrm{co}\hspace{0.2mm} H}$.
Notice from \eqref{ETPDop} and \eqref{eqRwedgeW} that we have an isomorphism
$\overline{H}\ot R\simeq H$ of $R$-algebras in ${}^{\overline{H}}\M$, which, with $D \square_{\overline{H}}$
applied, induces an isomorphism   
\begin{equation}\label{EqDRDH}
D \otimes R \simeq D\square_{\overline{H}}H
\end{equation}
of super-algebras over $B\otimes R$. We conclude that
$B\hookrightarrow D \square_{\overline{H}} H$, $b \mapsto b \otimes 1$ is faithfully flat. 

We claim that the inclusion $D \square_{\overline{H}} H\hookrightarrow D\ot H$ is faithfully flat.
To prove this, let $\psi : \overline{H}\to H$ be the restriction of an isomorphism $\overline{H}\ot R \simeq H$
such as above, to $\overline{H}\, (=\overline{H}\ot k)$. Then this $\psi$ is an algebra morphism in ${}^{\overline{H}}\M$.
Moreover, the isomorphism is given explicitly by $a \ot x \, (\in \overline{H}\ot R) \mapsto \psi(a)x$, and 
the induced one \eqref{EqDRDH} is given by
\[
d \otimes x \, (\in D \otimes R) \mapsto d_{(0)} \otimes \psi(d_{(1)})x,
\]
where $d \mapsto d_{(0)}\otimes d_{(1)}$ presents the $\overline{H}$-co-action on $D$.  
Compose \eqref{EqDRDH} with the inclusion in question. The resulting map is seen to be the composite of 
the inclusion $D \otimes R \hookrightarrow D \otimes H$ with the super-algebra map
\[
D \otimes H \to D \otimes H,\quad d \otimes h \mapsto d_{(0)}\otimes \psi(d_{(1)})h,
\]
which is indeed an isomorphism with inverse $d \otimes h \mapsto d_{(0)}\otimes \psi(\mathcal{S}(d_{(1)}))h$,
where $\mathcal{S}$ denotes the antipode of $\overline{H}$. 
The faithful flatness of the last inclusion proves the claim. 

It remains to prove that the beta map for $D \square_{\overline{H}} H$
is bijective. It suffices to show that the map turns to be so after the base extension along 
$D \square_{\overline{H}} H\hookrightarrow D\ot H$, which is faithfully flat as was just proven. 
To the beta map for $D$, apply first $(D\ot H)\ot_D$, and then $\square_{\overline{H}} H$. The resulting bijection
\[
(D\ot H)\ot_B(D \square_{\overline{H}}H) \overset{\simeq}{\longrightarrow} (D\ot H)\ot H
\]
is given explicitly by
\[
(d\ot h)\ot_B \big(\sum_i d_i \ot h_i\big)\mapsto \sum_i (-1)^{|d_i||h|} dd_i \ot h \ot h_i.
\]
This, composed with the bijection $D \ot H \ot H \overset{\simeq}{\longrightarrow} D \ot H \ot H$,\ $d \ot h \ot h' \mapsto d\ot h h'_{(1)}\ot h'_{(2)}$, 
coincides with the base extension which we wish to prove to be bijective; 
it is thus bijective, indeed.  
\epf

\begin{definition}\label{DARI}
An $H$-Galois extension $A$ over an arbitrary non-zero super-algebra $B$
is said to \emph{arise from} an $\overline{H}$-Galois extension $D/B$ if it is isomorphic
to the algebra $D \square_{\overline{H}} H$ in $\mathtt{SMod}^H$ over $B$. 
This $D \square_{\overline{H}} H$ is called the \emph{inflation} of $D$ along $H \to \overline{H}$. 
\end{definition}

For the $R$ in \eqref{ER}, let 
\[
R^+=R \cap \mathrm{Ker}(\varepsilon_H)
\]
denote the augmentation super-ideal, where 
$\varepsilon_H$ is the co-unit of $H$. 

\begin{prop}\label{PARI}
Let $B$ be a non-zero super-algebra. 
\begin{itemize}
\item[(1)] An $H$-Galois extension $A$ over $B$ arises from some $\overline{H}$-Galois
extension if and only if there exists an algebra morphism $R\to A$ in $\mathtt{SMod}^H$.
\item[(2)] Suppose that $A=(A, \rho)$ is an $H$-Galois extension over $B$, and
$f : R \to A$ is an algebra morphism in $\mathtt{SMod}^H$. Then the quotient super-algebra of $A$
\begin{equation}\label{EDF}
D_f:=A/(f(R^+))
\end{equation}
divided by the super-ideal $(f(R^+))$ generated by the image $f(R^+)$ of $R^+$
is a quotient algebra in $\mathtt{SMod}^{\overline{H}}$, and is in fact an $\overline{H}$-Galois extension over $B$. Moreover, 
\begin{equation}\label{EIOT}
\iota:A \longrightarrow D_f \, \square_{\overline{H}}H,\quad 
\iota(a)= \overline{a}_{(0)}\ot a_{(1)}
\end{equation}
is an isomorphism of $B$-algebras in $\mathtt{SMod}^H$, whence $A/B$ arises from $D_f/B$. 
Here, $\rho(a)=a_{(0)}\ot a_{(1)}$, and
$A \to D_f,\ a\mapsto \overline{a}$ denotes the quotient map.
\end{itemize}
\end{prop}

\pf
If $A=D\, \square_{\overline{H}}H$ is the inflation of $D$, then the map
\[ f : R \to A=D\, \square_{\overline{H}}H,\quad f(x)=1\ot x \]
is an algebra morphism in $\mathtt{SMod}^H$. 
This proves ``only if" of Part 1. 

Suppose that $A/B$ is $H$-Galois, and $f : R \to A$ is an algebra morphism in 
$\mathtt{SMod}^H$. Then $A$, regarded as an object in $\mathtt{SMod}_R$ through $f$,
turns into an object in $\mathtt{SMod}^H_R\, 
(=(\mathtt{SMod}^H)_R)$. 
Proposition 1.1 of \cite{M}
proves a category equivalence, 
\[
\M^{\overline{H}} \overset{\approx}{\longrightarrow} \M^H_R,\quad 
N\mapsto N\square_{\overline{H}}H,
\]
and that $\iota$ is an isomorphism. It remains to prove
that $D_f$ is $\overline{H}$-Galois over $B$. 
Now, $B \to D_f \, \square_{\overline{H}}H\, (=A)$ is faithfully flat. 
Notice from \eqref{ETPDop} that $N\square_{\overline{H}}H=N^{\oplus d}$ with $d=\dim(\wedge(W))$, whence
the functor $\square_{\overline{H}}H$
is faithfully exact.
Then it follows that $B \to D_f$ is faithfully flat. 
Apply $D_f\ot_A$ to the beta
map for $A$. 
Then there results an isomorphism in $\M^H_R$, which is seen to correspond,
through the category equivalence above, to the beta map for $D_f$;
it is, therefore, bijective.
\epf

\begin{rem}\label{RGUE}
One can generalize to our super situation, 
the result \cite[Theorem 4]{G} by G\"{u}nther proved in the non-super situation,
by modifying his proof. But the last proposition is formulated and proved 
under the restricted, super-commutativity assumption; it
allows us to ignore the Miyashita-Ulbrich actions used in \cite{G}, see
Page 4391, line 3.
\end{rem}

\subsection{Main theorems re-formulated in the Hopf-algebra language}\label{SSMR}
In this subsection we 
let $H$ be a Hopf super-algebra
which is assumed to be finitely generated unless otherwise stated.
We let $\G$ denote the affine super-group represented by $H$; it will be discussed only when $H$
is assumed to be finitely generated, and will be, therefore, an affine algebraic super-group.
Notice that the Hopf algebra $\overline{H}$ is finitely generated under the assumption.

\begin{prop}\label{PHA}
Let $\g=\Lie(\G)$ be the Lie super-algebra of $\G$, and
assume that odd component $\g_1$ of $\g$ satisfies
\begin{equation}\label{Eggg}
[[\mathfrak{g}_1,\mathfrak{g}_1],\mathfrak{g}_1]=0. 
\end{equation}
Suppose that $B$ is an arbitrary non-zero super-algebra. Then
every $H$-Galois extension over $B$ arises from an $\overline{H}$-Galois extension.
\end{prop}

The assumption above is clearly satisfied in case $[\g_1,\g_1]=0$, in which case the proposition 
will be applied to prove Lemma \ref{LHB} contained in Step 2 of the proof of
Theorem \ref{THB}; the theorem is a key to the proof of the second main theorem, Theorem \ref{THMII}.  
The assumption is also satisfied if $\dim \g_1=1$, since in general, a Lie super-algebra $\g$
shall satisfy $[[x,x],x]=0$ for all $x \in \g_1$. 
Proposition \ref{PHA} will be proved by spending the whole of Section \ref{SPHA}. 

Notice that the next theorem has a stronger conclusion than Proposition \ref{PHA} above;
see the sentence following the theorem.  

\begin{theorem}\label{THA}
Let $B\ne 0$ be a super-algebra. Assume
\begin{itemize}
\item[(i)] $\overline{H}$ is smooth as an algebra,\
\item[(ii)] the natural algebra surjection $B_0\to \overline{B}$ splits, and 
\item[(iii)] the super-ideal $(B_1)$ of $B$ is nilpotent.
\end{itemize}
Choose arbitrarily a section $\overline{B}\to B_0$ of the surjection $B_0\to \overline{B}$, and regard $B$ as a super-algebra
over $\overline{B}$ through it. Then every $H$-Galois extension $A$ over $B$ 
is isomorphic to 
\begin{equation}\label{EBAH}
(B\otimes_{\overline{B}}\overline{A})\square_{\overline{H}} H.
\end{equation}
Moreover, $\overline{A}$, which is regarded naturally as a purely even algebra in 
$\M^{\overline{H}}$, is a purely even $\overline{H}$-Galois extension over $\overline{B}$. 
\end{theorem}

As a result we can say that $A/B$ arises from $B\otimes_{\overline{B}}\overline{A}$,
which is clearly an $\overline{H}$-Galois extension 
over $B$. 

\begin{rem}\label{RSmooth}
(1)\
As for the assumption (i) above, it is proved by \cite[Proposition A.3]{MZ2}
without assuming $H$ is finitely generated that
the following are equivalent:
\begin{itemize}
\item[(a)]
$\overline{H}$ is smooth as an algebra;
\item[(b)]
$H$ is smooth as a super-algebra.
\end{itemize}

These equivalent conditions are necessarily satisfied if $\operatorname{char} k=0$. Indeed, 
it is well known that $\overline{H}$ (in $\operatorname{char} k=0$) is geometrically reduced; this
is equivalent to saying that $\overline{H}$ is smooth
in the restricted case where it is finitely generated. In the general case, $\overline{H}$ is 
a filtered union of finitely generated Hopf subalgebras which are smooth as was just seen, and it is, therefore, 
smooth; see \cite[Exercise 9.3.2, 4]{Wei}.

(2)\
Example 3.23 of \cite{MSS} shows that 
Theorem \ref{THA} does not hold without the assumption (i), even if the remaining 
(ii) and (iii) are satisfied. 

(3)\
$(B\otimes_{\overline{B}}\overline{A})\square_{\overline{H}} H$ in \eqref{EBAH}
is canonically isomorphic to
$B\otimes_{\overline{B}}(\overline{A}\, \square_{\overline{H}} H)$, as is easily seen by using
\eqref{ETPDop}.
\end{rem}

Theorem \ref{THA} will be proved by spending the whole of
Section \ref{STHA}. The theorem has the following corollary,
which is seen to be translated into the first main theorem, Theorem \ref{TBij}, 
in view of Corollary \ref{CXGAB} and Remark \ref{RSmooth} (1) just above. 

\begin{corollary}\label{CHA}
Assume that $\overline{H}$ is smooth. 
Let $B\ne 0$ be a Noetherian smooth algebra. If $A/B$ is $H$-Galois, then 
$\overline{A}/\overline{B}$ is purely even $\overline{H}$-Galois. 
Moreover, the assignment $A \mapsto \overline{A}$ gives rise to a bijection
from
\begin{itemize}
\item the set of all isomorphism classes of $H$-Galois extensions over $B$
\end{itemize}
onto
\begin{itemize}
\item the set of all isomorphism classes of purely even $\overline{H}$-Galois extensions 
over $\overline{B}$.
\end{itemize}
\end{corollary}
\pf
To a purely even $\overline{H}$-Galois extension $E/\overline{B}$, one can assign 
an $H$-Galois extension, $(B\otimes_{\overline{B}}E)\square_{\overline{H}}H$, over $B$.
This assignment gives rise to 
a desired inverse of $A\mapsto \overline{A}$,
as is easily seen by using \eqref{ETPDop} and Theorem \ref{THA}.
\epf

\begin{rem}\label{RNoncat}
We have an example of an $H$-Galois extension $(B\otimes_{\overline{B}}E)\square_{\overline{H}}H$ such as in the above proof, with the property that
it has a $B$-algebra endomorphism
(in fact, an automorphism by Lemma \ref{LIsom}) in $\M^H$ which is not induced naturally 
from any $\overline{H}$-comodule $\overline{B}$-algebra endomorphism of $E$. Indeed, let
$H=\overline{H}=k[t,t^{-1}]$, the Laurent polynomial algebra regarded as the Hopf algebra
which represents the multiplicative group. Let $B=\wedge(V)$ be the exterior algebra
on a finite-dimensional vector space of dimension $>1$. Choose as $E$ the trivial
$k[t,t^{-1}]$-Galois extension $k[t,t^{-1}]$
over $k$. Then one has $(B\otimes_{\overline{B}}E)\square_{\overline{H}}H=\wedge(V)\otimes k[t,t^{-1}]$. Every endomorphism in question
uniquely arises from a super-algebra map $g : k[t,t^{-1}] \to \wedge(V)$ so that
\[
\wedge(V)\otimes k[t,t^{-1}]\to \wedge(V)\otimes k[t,t^{-1}],\quad b \otimes h \mapsto 
b(g \otimes \mathrm{id})(\Delta(h)).  
\]
The endomorphism that arises from the $g$ defined by $g(t^{\pm 1})=1\pm x$, where $0\ne x \in \wedge^2(V)$ such that $x^2=0$, is not induced from 
any endomorphism of $E=k[t,t^{-1}]$, as is easily seen. 
The example shows that the bijection proved
in Corollary \ref{CHA} does not result from a category-equivalence; this is translated 
into the fact stated in Remark \ref{RNoHope} (1).
\end{rem}

The second main theorem, Theorem \ref{THMII}, restricted to the
affinity situation, is translated into the following.

\begin{theorem}\label{THB}
Suppose that $A=(A, \rho)$ is a non-zero algebra in $\M^H$; 
$\overline{A}$ is then naturally a purely even algebra in $\M^{\overline{H}}$.  
Let $B=A^{\mathrm{co}\hspace{0.2mm}H}$. 
Assume
\begin{itemize}
\item[(a)] 
the map $\alpha : A\otimes A\to A\otimes H$ in \eqref{EALP} is surjective,
\item[(b)] 
$\overline{A}/\overline{A}^{\mathrm{co}\hspace{0.2mm}\overline{H}}$ is $\overline{H}$-Galois,
\item[(c)] 
$\overline{H}$ is smooth as an algebra, and 
\item[(d)] 
$A$ is Noetherian and smooth.
\end{itemize}
Then $A/B$ is $H$-Galois, and $B$ is Noetherian, smooth and such that
\begin{equation}\label{EqBAH}
\overline{B}=\overline{A}^{\mathrm{co}\hspace{0.2mm}\overline{H}},
\end{equation}
whence, with the assumptions of Theorem \ref{THA} all satisfied, $A$ is of the form
\eqref{EBAH}.  
\end{theorem}

\begin{rem}\label{RHB}
Under the situation of Theorem \ref{THB} assume that (a) is satisfied. It is
proved by \cite[Theorems 1.2 and 3.7]{MSS} that $A/B$ is $H$-Galois, if $\overline{H}$ 
is co-Frobenius as a coalgebra,
or equivalently, if the affine algebraic group
$\G_{\mathsf{ev}}$, which is now represented by $\overline{H}$,
has an integral. 
When $\operatorname{char}k=0$, the assumption 
is satisfied if and only if $\G_{\mathsf{ev}}$ is linearly reductive;
see Sullivan \cite{Sullivan} also for
characterizations in the positive characteristic case. 
\end{rem}

Theorem \ref{THB} will be proved in Section \ref{STHB}, and the result, regarded 
as the special, affinity case of Theorem \ref{THMII}, will be applied to
prove that theorem in Section \ref{STHMII}. 

In what follows until the end of this subsection we let $H$ be a Hopf super-algebra
which may not be finitely generated.

The notion of smoothness for super-algebras (over $k$) is generalized to the notion for 
super-algebras over a fixed super-algebra
and, moreover, for those super-algebras on which 
$H$-co-acts. 
For the latter we 
use the word ``$H$-smooth", as follows.

\begin{definition}\label{DHSmoo}
Let us be given an algebra morphism $R \to S$ in $\M^H$. We say that $S$ is 
$H$-\emph{smooth over} $R$, if given a surjective
$R$-algebra morphism $A \to B$ in $\M^H$ with nilpotent kernel, every $R$-algebra
morphism 
$S \to B$ in $\M^H$ factors through $A$ by some $R$-algebra morphism $S \to A$ in $\M^H$. 
The condition is equivalent to saying that every surjective $R$-algebra morphism in $\M^H$ that maps
onto $S$ and has nilpotent kernel necessarily splits, as is seen by a familiar augment using pull-backs; 
see the proof of \cite[Proposition 9.3.3]{Wei}, for example. Here one should notice 
that the category of $R$-algebras in $\M^H$ has pull-backs.
\end{definition}

We will use the word, when $R$ such as above is trivial as an $H$-super-comodule 
(and, moreover, for applications in the sequel, when $R \to S$ is an inclusion and $R\subset S^{\mathrm{co}\hspace{0.2mm}H}$).

\begin{prop}\label{PHSmSm}
Let us be given an algebra morphism $R \to S$ in $\M^H$ such that
$R$ is trivial as an $H$-super-comodule. If $S$ is $H$-smooth over $R$,
then it is necessarily smooth over $R$.
\end{prop}
\pf
Let $T$ be an $R$-super-algebra in general. Then $T \otimes H$ is 
naturally an $R$-algebra in $\M^H$. By the triviality assumption for $R$, 
the $R$-super-algebra maps $S \to T$ are in one-to-one correspondence 
with the $R$-algebra morphisms $S \to T \otimes H$ in $\M^H$, in which
an $R$-super-algebra map $f : S \to T$ corresponds to 
\begin{equation}\label{AddedEq}
\widetilde{f}:=(f\otimes \operatorname{id}_H)\circ \rho : 
S \to T \otimes H,
\end{equation}
where $\rho: S\to S \otimes H$ denotes the $H$-co-action on $S$.

Given an $R$-super-algebra surjection $A \to B$ with nilpotent kernel, 
we have a surjective $R$-algebra morphism $A\otimes H \to B\otimes H$ in $\M^H$
with nilpotent kernel. Assume that $S$ is $H$-smooth over $R$. Given 
an $R$-super-algebra map $f : S \to B$, the $R$-algebra morphism 
$\widetilde{f} : S \to B\otimes H$ in $\M^H$ factors through $A \otimes H$
by some $S \to A \otimes H$ that is necessarily of the form 
$\widetilde{g}$ with $g : S \to A$ an $R$-super-algebra map. 
We see that $f$ factors through $A$ 
by $g$. This proves the desired smoothness. 
\epf

By the proposition just proven, $H$ is smooth (over $k$), if it is 
$H$-smooth over $k$. The converse holds, as well, as is shown by the following.

\begin{prop}\label{PHSmooth}
$H$ \textup{(}or equivalently, $\overline{H}$\textup{)} 
is smooth if and only if it is $H$-smooth over $k$. 
These equivalent conditions are necessarily 
satisfied if $\operatorname{char} k=0$; see Remark \ref{RSmooth} (1).
\end{prop}
\pf
By modifying in an obvious manner the proof of the result 
\cite[Proposition 1.10, (i)$\Leftrightarrow$(iii)]{MO}
in the non-super situation, we see that the following are equivalent:
\begin{itemize}
\item[(a)] For every trivial $H$-super-module $M$, the symmetric 2nd Hochschild cohomology $H_s^2(H,M)$
(constructed in $\M_k$) vanishes; 
\item[(b)] Every surjective algebra morphism in $\M^H$ that maps onto $H$ and
has nilpotent kernel necessarily splits.
\end{itemize}

We now have only to prove ``only if".  If $H$ is smooth, then (a) and thus (b) are satisfied.
But (b) is equivalent to saying that $H$ is $H$-smooth over $k$,
as was remarked at the last part of of Definition \ref{DHSmoo}.
\epf

\begin{rem}\label{AddedRem}
Let $H'$ be a quotient Hopf algebra of $H$ such that $H$ is injective as a left or equivalently, right
$H'$-super-comodule; recall from the second paragraph of Section \ref{SSSCA} that 
the condition is equivalent to saying that
the co-tensor product functor $\square_{H'} H$ or $H \square_{H'}$ is exact. In the situation of 
Proposition \ref{PHSmSm}, if $S$ is $H$-smooth over $R$, then $S$ is $H'$-smooth over $R$.
This improves Proposition \ref{PHSmSm}, generalizing $k\, (=H/\operatorname{Ker}\varepsilon)$ to $H'$.
But the proof is a slight modification, which replaces \eqref{AddedEq} with
$(f \square_{H'} \operatorname{id}_H)\circ \rho : S \to T \square_{H'} H$.

It follows that the conditions proved by Proposition \ref{PHSmooth} to be equivalent to each other 
are further equivalent to the condition that $H$ is $H'$-smooth for some/any quotient Hopf super-algebra
$H'$ of $H$ that satisfies the injectivity assumption.
\end{rem}

\section{Proof of Proposition \ref{PHA}}\label{SPHA}

To prove the proposition (in 2 steps), let $H$ be a finitely generated Hopf super-algebra,
let $\G$ denote affine algebraic super-group represented by $H$, and let
$\g=\Lie(\G)$ be the Lie super-algebra of $\G$. We assume $[[\g_1,\g_1],\g_1]=0$
as in \eqref{Eggg}. 

\subsection{Step 1}\label{SPA1}
We claim the following.

\begin{lemma}\label{LGF}
$\G$ has a quotient affine super-group $\mathsf{F}$ which has the following two properties, 
where we let 
$\mathfrak{f}=\Lie(\mathsf{F})$ denote the Lie super-algebra of $\mathsf{F}$.
\begin{itemize}
\item[(1)] The natural Lie super-algebra map $\g \to \mathfrak{f}$, restricted to
the odd components, is identical, or in notation, $\g_1=\mathfrak{f}_1$;
\item[(2)] The closed embedding $\mathsf{F}_{\mathsf{ev}} \hookrightarrow \mathsf{F}$ splits
\textup{(}necessarily, uniquely\textup{)}. 
\end{itemize}
\end{lemma}

Suppose that this lemma is proved. Let $J=\cO(\mathsf{F})$. Then $J$ is a Hopf
super-subalgebra of $H$. By (2), this $J$ is of the form of smash co-product
\begin{equation}\label{EJJW}
J = \overline{J}~\rcosmash\, \wedge(W).
\end{equation}
Here $W$ is a purely odd object in $\M^{\overline{J}}$; it
constitutes the Hopf algebra $\wedge(W)$ in $\M^{\overline{J}}$, which
contains elements of $W$ as primitives. The Hopf algebra in turn
constitutes
the smash co-product above.
By (1), we have
$W=(\mathfrak{f}_1)^*= (\g_1)^*$, whence
\[ {}^{\mathrm{co}\hspace{0.2mm}\overline{H}}H=\wedge(W). \]
Note that the $k\oplus W$ in this ${}^{\mathrm{co}\hspace{0.2mm}\overline{H}}H$ 
is a sub-object of $H$ in $\M^H$, since it is, in fact, such of $J$ in $\M^J$. 

Given an $H$-Galois extension $A/B$, the inclusion $B \hookrightarrow A$
has a retraction $r : A \to B$ in $\M_B$. Indeed, by applying
\cite[Theorem 4.9]{SS} to the $\widehat{H}$-Galois extension $\widehat{A}/B$
in the non-super situation which was discussed in the proof of Proposition \ref{PSS}, 
one sees
that the natural inclusion $B \hookrightarrow \widehat{A}\, (=A\rtimes \mathbb{Z}_2)$
splits left (or right) $B$-linearly. The argument using separable
extensions, which was used in the last-mentioned proof, ensures the existence of $r$. 
Now, the composite
\[
k\oplus W \hookrightarrow H=k\ot H
\hookrightarrow A\ot H \overset{\beta^{-1}}{\longrightarrow}A\ot_BA 
\overset{r\ot_B\op{id}_A}{\longrightarrow} A
\]
is a morphism in $\M^H$ which sends $1$ in $k\, (\subset k\oplus W)$  to $1 \in A$. 
Therefore, this uniquely extends to an algebra morphism 
${}^{\mathrm{co}\hspace{0.2mm}\overline{H}}H=\wedge(W)\to A$ in $\M^H$. 
Proposition \ref{PARI} (1) then proves Proposition \ref{PHA} in question.

\subsection{Step 2}\label{SPA2}
It remains to prove Lemma \ref{LGF} above. Let
\[
C=H^{\circ},\quad \underline{C}=\overline{H}^{\circ}
\]
denote the dual Hopf super-algebra of $H$ and the dual Hopf algebra of $\overline{H}$;
see \cite[Sect.~2.4]{M1}. 
The Hopf super-algebra $C$ is not necessarily super-commutative, but is super-cocommutative.
This $C$ includes $\underline{C}$ as the largest purely even Hopf super-subalgebra. The spaces
of all primitives in the Hopf (super-)algebras give the Lie (super-)algebras of the
affine (super-)groups, so as
\[
P(C)=\Lie(\G)\, (=\g),\quad P(\underline{C})=\Lie(\G_{\mathsf{ev}})\, (=\g_0).
\]

Let $\langle \ , \ \rangle : \underline{C}\times \overline{H} \to k$ denote the canonical pairing. Then $\overline{H}$ has the natural left $\underline{C}$-module
structure given by
\[
c\rightharpoonup a = a_{(1)}\langle c, a_{(2)}\rangle,
\]
where $c \in \underline{C}$,\ $a \in \overline{H}$, and $\Delta(a)=a_{(1)}\ot a_{(2)}$
denotes the co-product on $\overline{H}$. 
It is shown by \cite[Lemma 27, Theorem 29]{M1} that the set
\begin{equation}\label{EHom}
\op{Hom}_{\underline{C}}(C,\overline{H})
\end{equation}
of all left $\underline{C}$-linear maps $C\to \overline{H}$ has a natural 
structure of Hopf super-algebra, such that 
$H \simeq \op{Hom}_{\underline{C}}(C,\overline{H})$ naturally as Hopf super-algebras.
This Hom set is essentially the same as what is denoted by 
$\op{Hom}_J(H(J,V), C)$ in \cite[Page 1101, line 11]{M1}, and the 
explicit Hopf-super-algebra structure is given by
\cite[Proposition 18 (2), (3)]{M1}. 

Let
\[ V = [\g_1,\g_1]\, (\subset \underline{C}). \]
One sees that under the adjoint actions, 
$\underline{C}$ stabilizes $V$, and $\g_1$ annihilates $V$ by 
\eqref{Eggg}. It follows that the left and the right (super-)ideals of the
Hopf (super-)algebra generated by $V$ coincide, so as
\[
\underline{C}V=V\underline{C},\quad CV=VC.
\]
Let
\[
\underline{Q}=\underline{C}/V\underline{C},\quad Q=C/VC
\] 
denote the resulting quotients. 
In view of \cite[Theorem 3.6]{M} one sees that $Q$
is a super-cocommutative (but not necessarily,
super-commutative)
Hopf super-algebra, which includes
$\underline{Q}$
as the largest purely even Hopf super-subalgebra,
and the odd component of the
primitives $P(Q)$ in $Q$ remains to be $\g_1$, or in notation,
\begin{equation}\label{EPQ1}
P(Q)_1=\g_1.
\end{equation}
Moreover, $[\g_1,\g_1]=0$ in $P(Q)$, whence 
$Q$ is of the form of smash product
\[
Q=\underline{Q}~\ltimes \, \wedge(\g_1). 
\]

Define
\[
\overline{H}':=\{\, a \in \overline{H}\, \mid \, V\rightharpoonup a=0 \, \}.
\]
Then this is naturally a left $\underline{Q}$-module, whence we have the Hom-set
\begin{equation}\label{EHom1}
J:=\op{Hom}_{\underline{Q}}(Q,\overline{H}')
\end{equation}
as a subset of the one in \eqref{EHom}. 
Since $\overline{H}$
is finitely generated, it is proper in the sense that the algebra map
$\overline{H}\to (\underline{C})^*$ which arises from the canonical
pairing is injective. Therefore, we have
\[
\overline{H}'=\{\, a \in \overline{H} \, \mid \, \langle \underline{C}V,\ a \rangle =0\, \}.
\]
Moreover, since $\underline{C}V$ is a Hopf ideal of $\underline{C}$, 
it follows that $\overline{H}'$
is a Hopf subalgebra of $\overline{H}$. 
From the above-mentioned structure on 
$\op{Hom}_{\underline{C}}(C,\overline{H})$,
we see that the $J$ in \eqref{EHom1} is regarded as a Hopf super-subalgebra of 
$H$, such that
\[
(\overline{J}=)\, J/(J_1)=\op{Hom}_{\underline{Q}}(\underline{Q}\ot k, \overline{H}') =\overline{H}'.
\]
In addition, the natural embedding
\[
\overline{J}=\op{Hom}_{\underline{Q}}(Q/Q \g_1, \overline{H}')\ \hookrightarrow \
\op{Hom}_{\underline{Q}}(Q,\overline{H}')
\]
is seen to be a Hopf super-algebra section of the projection $J \to \overline{J}$. 
It follows that $J$ represents a quotient affine super-group, say $\mathsf{F}$, of $\G$
which has the property (2) of Lemma \ref{LGF}. 
We see from \eqref{EPQ1} that it has the property (1), as well.


\section{Proof of Theorem \ref{THA}}\label{STHA}

To prove the theorem, let $H$ be a Hopf
super-algebra, 
and suppose that $A/B$ is an $H$-Galois extension. 
The proof is divided into 3 steps. In Steps 2 and 3 we will assume that $H$ is finitely generated.

We start with only assuming (i), or namely, that 
$\overline{H}$ is smooth over $k$. 

\subsection{Step~1}\label{STHA1}
This is a crucial step devoted to proving 
the following proposition, for which, we emphasize, $H$ may not be finitely generated.

\begin{prop}\label{PHS}
$A$ is $H$-smooth over $B$, whence it is smooth over $B$ by Proposition \ref{PHSmSm}. 
\end{prop}

We let $\A$ denote the super-algebra $A$ on which $H$ co-acts trivially; thus we have
$(\A)^{\mathrm{co}\hspace{0.2mm} H}= \A$. 

By Proposition \ref{PHSmooth}, $H$ is $H$-smooth. This implies the following.

\begin{lemma}\label{LHS}
$\A\ot H$ is $H$-smooth over $\A$. 
\end{lemma}

Define
\[
\bA:=\A\ot_B A.
\]
Note that $H$ co-acts on the second tensor factor $A$. Since this is isomorphic
to $\A \ot H$ through the beta map in \eqref{EBET}, the lemma above has the following corollary.

\begin{corollary}\label{CHS}
$\bA$ is $H$-smooth over $\A$.
\end{corollary}

Recall from \eqref{ESMod} the construction of the category 
$\M^H_A$, which is in fact $k$-linear abelian. 
By an analogous construction we let 
\begin{equation}\label{EAMA}
{}_A^{}\underset{B}{\M}{}_A^H
\end{equation}
denote
the $k$-linear abelian category of those $(A,A)$-bimodules $M$ in $\M^H$ for which the left and the right
$B$-actions coincide (in the super sense), or explicitly,
\[ bm=(-1)^{|b||m|}mb,\quad b\in B,\ m\in M. \]
For example, $A\ot_B A$, equipped with the co-diagonal $H$-co-action, is an object
of this category. Clearly, $A$ is, as well. For simplicity we will write 
${}_A^{}\M{}_A^H$ for \eqref{EAMA}. Notice that the category includes 
$\M^H_A$ as a full subcategory; it consists of all
objects of ${}_A^{}\M{}_A^H$ for which the left and the right $A$-actions coincide.

Recall from Theorem \ref{TOBER} that we have the ($k$-linear) category-equivalence
\begin{equation}\label{EPHI}
\Phi : \M_A^H \to \M_B,\quad \Phi(M)=M^{\mathrm{co}\hspace{0.2mm} H}.
\end{equation}

\begin{lemma}\label{LEXT}
The direct sum
\[
\mathbb{P}:=(A\otimes_BA)[0]\oplus (A\otimes_BA)[1]
\]
of $A\otimes_BA\, (=(A\otimes_BA)[0])$ and its degree shift $(A\otimes_BA)[1]$
is a projective generator in ${}_A^{}\M{}_A^H$, whence 
the category has enough projectives, and for every $M \in {}_A^{}\M{}_A^H$, the Ext group 
\[ \mathrm{Ext}^\bullet_{{}_A^{}\M{}_A^H}(-, M)  \]
is defined as the right derived functor of the 
Hom functor $\mathrm{Hom}_{{}_A^{}\M{}_A^H}(-, M)$. 
\end{lemma}
\pf
Let $M \in {}_A^{}\M{}_A^H$. We have the natural identifications
\[
\mathrm{Hom}_{{}_A^{}\M{}_A^H}((A\ot_BA)[i], M)= \big(M^{\mathrm{co}\hspace{0.2mm} H}\big)_i, \quad i=0,1
\]
of $k$-vector spaces, 
which amount to
\begin{equation}\label{EqPM}
\mathrm{Hom}_{{}_A^{}\M{}_A^H}(\mathbb{P}, M)= M^{\mathrm{co}\hspace{0.2mm} H}.
\end{equation}
With $M$ regarded as an object in $\M{}_A^H$
through the obvious forgetful functor, the right-hand side of \eqref{EqPM}
coincides with $\Phi(M)$; here, we do not regard 
$\M{}_A^H\subset {}_A\M{}_A^H$.
Since $\Phi$ preserves surjections, it follows that $\mathbb{P}$ as well as
the two direct summands are projective. 
Since the inclusion $\Phi(M)=M^{\mathrm{co}\hspace{0.2mm} H}\hookrightarrow M$
extends to a natural isomorphism $M^{\mathrm{co}\hspace{0.2mm} H}\otimes_BA\overset{\simeq}{\longrightarrow} M$ in $\M{}_A^H$, one sees that $M$ is generated
by $M^{\mathrm{co}\hspace{0.2mm} H}$ as a right $A$-module, and so as an $(A,A)$-bimodule.
Therefore, $\mathbb{P}$ is a generator. 
\epf

\begin{rem}\label{REXT}
Since $\bA=\A\otimes_B A$ is an $H$-Galois extension over $\A$, an analogous result of the
lemma above holds, with $A\otimes_B A$ replaced by
$\bA~{\otimes}_{\A}\bA\, (=\A \otimes_B(A\otimes_B A))$, in
the $k$-linear abelian category
\[
{}_\bA^{}\underset{\A}{\M}{}_\bA^H,
\]
which is defined in an obvious manner. 
\end{rem}

Note that $B$ is a subalgebra of $A$ in $\M^H$. 
Just as shown in \cite[Section 9.3]{Wei} in the ordinary situation 
(see also the last part of Definition \ref{DHSmoo}),
$A$ is $H$-smooth over $B$ if and only if 
an \emph{extension of commutative $B$-algebras in} $\M^H$
\begin{equation}\label{EMPA}
0 \to M \to P \to A \to 0 
\end{equation}
necessarily splits for every $M \in \M{}_A^H\, (\subset {}_A^{}\M{}_A^H)$. By an \emph{extension} 
such as above we mean that $P \to A$ is a 
surjective morphism of
commutative $B$-algebras 
in $\M^H$ with kernel $M$ square-zero, $M^2=0$, in $P$. It is said to \emph{split} if
the last surjection splits as a $B$-algebra morphism in $\M^H$.  

\begin{lemma}\label{LSP}
Every $B$-algebra extension of $A$, such as \eqref{EMPA}, necessarily splits as 
a short exact sequence of $B$-modules in $\M^H$. 
\end{lemma}
\pf
By Proposition \ref{PSS} the product map $B \otimes A \to A$ splits
as a $B$-module morphism in $\M^H$. 

It follows that every object in $\M_A^H$, and thus the kernel $M$, in particular, are injective in $\M^H$.
Indeed, we see that $M=M^{\mathrm{co}\hspace{0.2mm} H}\otimes_B A$ is 
a direct summand of 
$M^{\mathrm{co}\hspace{0.2mm}H}\otimes_B(B\otimes A)
=M^{\mathrm{co}\hspace{0.2mm} H}\otimes A$ in $\M^H$.
The desired injectivity follows since $A$
is injective in $\M^H$ by
Theorem \ref{TOBER} (1). 

We thus have a section $s : A \to P$ in $\M^H$
of the surjection $P \to A$.
Choose arbitrarily a $B$-module 
section $t : A \to B \otimes A$ in $\M^H$ of the product map above. 
We see that the composite
\[ A \overset{t}{\longrightarrow} B\otimes A\overset{\mathrm{id}\otimes s}{\longrightarrow} B\otimes P \to P, \]
where the last arrow indicates the product map, is a 
section of the $B$-module morphism $P \to A$ in $\M^H$. 
\epf

An important consequence of the lemma above is that
$B$-algebra extensions of $A$ such as above are classified
by the symmetric 2nd Hochschild cohomology $H_s^2(A,M)$ 
{(constructed in ${}_B\M^H$),
which is naturally embedded into
\[
\mathrm{Ext}^2_{{}_A^{}\M{}_A^H}(A, M).
\]
Let 
\[
J_{A/B}=\operatorname{Ker}(A\ot_BA\to A)
\]
be the kernel of the product map. Then one has the short exact sequence 
$0 \to J_{A/B}\to A\otimes_B A \to A \to 0$ in ${}_A^{}\M{}_A^H$. 
The derived long exact sequence, combined with the projectivity result in Lemma \ref{LEXT},
shows that the last 2nd Ext group is naturally isomorphic to
\[
\mathrm{Ext}^1_{{}_A^{}\M{}_A^H}(J_{A/B}, M). 
\]
Therefore, the algebra extension \eqref{EMPA} splits 
if and only if the corresponding
extension 
\begin{equation}\label{EQJ}
0\to M \to Q \to J_{A/B} \to 0 
\end{equation}
in ${}_A^{}\M{}_A^H$ splits;
in other words, the cohomology class of \eqref{EMPA} in $H_s^2(A,M)$
is zero if and only if its image in $\mathrm{Ext}^1_{{}_A^{}\M{}_A^H}(J_{A/B}, M)$
through the injection
\[
H_s^2(A,M)\hookrightarrow \mathrm{Ext}^2_{{}_A^{}\M{}_A^H}(A, M)
\overset{\simeq}{\longrightarrow} \mathrm{Ext}^1_{{}_A^{}\M{}_A^H}(J_{A/B}, M)
\]
is zero. 
By the flat base-extension $\A\otimes_B$, we obtain the mutually
equivalent, analogous conditions, and see that they are both satisfied since
$\bA = \A\otimes_BA$ is $H$-smooth over $\A$ by Corollary \ref{CHS}. 
The second analogous
condition tells us that the base extension
\begin{equation}\label{EAQJ}
0\to\A\otimes_B M \to \A\otimes_B Q \to \A\otimes_BJ_{A/B}=J_{\bA/\A} \to 0 
\end{equation}
of \eqref{EQJ} splits in ${}_{\bA}^{}\underset{\A}{\M}{}_{\bA}^H$.

Now, we let $\A$ recover the original structure as an algebra in $\M^H$, and suppose
that the sequence \eqref{EAQJ} has an additional structure which arises from
the tensor factors $\A$ regarded as objects in 
\[ {}_A^{}\M^H\, (={}_A^{}(\M^H)). \]
To make this
clearer, let $A_i$ (resp., $H_i$), $i=1,2$, be copies of $A$ (resp., $H$), 
regard $A_i$ as an algebra in $\M^{H_i}$, and suppose that $A_1$ is the $\A$
equipped with the recovered structure. One sees that
\[ A_1\ot_BA_2\quad \text{and}\quad A_2\, (=B\ot_BA_2) \]
naturally form algebras in $\M^{H_1\ot H_2}$, which both include $B$ as subalgebra.  
Therefore, we have the $k$-linear
abelian category
\[
{}_{A_1\ot_BA_2}^{}\underset{B}{\M}{}_{A_2}^{H_1\ot H_2}
\]
of those $(A_1\ot_BA_2,A_2)$-bimodules in $\M^{H_1\ot H_2}$ for which the left and the
right $B$-actions coincide. 
The short exact sequence \eqref{EAQJ} is naturally regarded as such a sequence in
${}_{A_1\ot_BA_2}^{}\underset{B}{\M}{}_{A_2}^{H_1\ot H_2}$. 

\begin{lemma}\label{LSPR}
The short exact sequence \eqref{EAQJ} splits 
in ${}_{A_1\ot_BA_2}^{}\underset{B}{\M}{}_{A_2}^{H_1\ot H_2}$. 
\end{lemma}

\pf
Our argument below is essentially the same as the one which proves 
the Maschke-type Theorem \cite[Theorem 1]{D}. 

Let us consider the injection $A_1\ot_BM\to A_1\ot_B Q$ in \eqref{EAQJ}.
Since the category ${}_{\bA}^{}\underset{\A}{\M}{}_{\bA}^H$ is naturally identified
with ${}_{A_1\ot_BA_2}^{}\underset{B}{\M}{}_{A_2}^{H_2}$, 
the injection has a retraction, say $r$, in
${}_{A_1\ot_BA_2}^{}\underset{B}{\M}{}_{A_2}^{H_2}$. The tensor product 
\[
\widetilde{r}:=r\ot \mathrm{id}_{H_1}:(A_1\ot H_1)\ot_B Q \to (A_1\ot H_1)\ot_BM.
\]
is seen to be in ${}_{A_1\ot_BA_2}^{}\underset{B}{\M}{}_{A_2}^{H_1\ot H_2}$,
where $A_1$ acts on the $A_1\ot H_1$ in these objects through the $H_1$-co-action on $A_1$,
and $H_1$ co-acts on the tensor factor $H_1$; to be precise, the latter co-action 
involves the super-symmetry relevant to $Q$ or $M$. 
One sees that $\tilde{r}$ is a retraction, 
in ${}_{A_1\ot_BA_2}^{}\underset{B}{\M}{}_{A_2}^{H_1\ot H_2}$, of the above injection 
tensored with $\mathrm{id}_{H_1}$. 

By Theorem \ref{TOBER} (1), $A_1$ is injective in $\M^{H_1}$. It follows that
there is a morphism $\phi : H_1 \to A_1$ in $\M^{H_1}$ such that $\phi(1)=1$. 
We see that the map
$t :A_1\ot H_1 \to A_1$ defined by
\begin{equation}\label{ERET}
t(a\ot h)=a_{(0)}\, \phi(\mathcal{S}(a_{(1)})h),
\end{equation}
where $\mathcal{S}$ denotes the antipode of $H_1$, 
is a retraction of the $H_1$-co-action $A_1 \to A_1\ot H_1,\ 
a\mapsto a_{(0)}\ot a_{(1)}$ in 
${}_{A_1}\M^{H_1}$; cf. \cite[Page 100, line --1]{D}. 
The tensor product
\[
\widetilde{t}:=t\ot_B \mathrm{id}_M: (A_1\ot H_1)\ot_BM\to A_1 \ot_B M
\]
is a retraction 
in ${}_{A_1\ot_BA_2}^{}\underset{B}{\M}{}_{A_2}^{H_1\ot H_2}$
of the $H_1$-co-action on $A_1\ot_B M$.
One sees that the composite $\tilde{t}\circ \tilde{r}$, further composed with
the $H_1$-co-action on $A_1\ot_B Q$, gives a retraction in ${}_{A_1\ot_BA_2}^{}\underset{B}{\M}{}_{A_2}^{H_1\ot H_2}$ of the injection in question. This proves the lemma. 
\epf

Apply to \eqref{EAQJ} the category equivalence
$(\ )^{\mathrm{co}(H_1)}: 
{}_{A_1}^{}\M{}^{H_1}\approx
{}_{B}\M$, which is essentially the same as the $\Phi$ in \eqref{EPHI}.
This can be, in fact, replaced by
\[
{}_{A_1\ot_BA_2}^{}\underset{B}{\M}{}_{A_2}^{H_1\ot H_2}\approx
{}_{A_2}^{}\M{}_{A_2}^{H_2}\, (={}_{A_2}^{}\underset{B}{\M}{}_{A_2}^{H_2}).
\]
Then there results
precisely the short exact sequence \eqref{EQJ} in the category 
${}_A^{}\M{}_A^H$ on the right-hand side, which
splits since by  Lemma \ref{LSPR},\ 
\eqref{EAQJ} splits in the category on the left-hand side.
This implies that $A$ is $H$-smooth over $B$, as was seen before. 

\subsection{Step~2}\label{STHA2}
Let us be in the situation at the beginning of this Section \ref{STHA}. 
But we here assume that $H$ is finitely generated (and continue to assume (i)).
In addition, we assume that $B$ is purely even, $B=B_0$, or equivalently,
$B =\overline{B}$, and aim to prove Theorem \ref{THA},
that is, the following. 

\begin{prop}\label{PPEB}
Under the assumptions above, we have
\[ A \simeq \overline{A}\, \square_{\overline{H}}H \]
\textup{(}or namely, $A$ is of the form \eqref{EBAH} with $B =\overline{B}$\textup{)}, 
and $\overline{A}/B$ is $\overline{H}$-Galois. 
\end{prop}

\begin{lemma}\label{LBF}
We have the following.
\begin{itemize}
\item[(1)]
$\overline{A}$ is smooth over $B$. 
\item[(2)]
$\overline{A}$ is faithfully flat over $B$.
\end{itemize}
\end{lemma}
\pf
(1)
By Proposition \ref{PHS}, $A$ is smooth over the now purely even $B$. 
As a general fact, this implies that 
$\overline{A}$ is smooth over $B$.
Indeed,
given an algebra surjection $S\to T$ with nilpotent kernel, every algebra
map $\overline{A}\to T$,
identified with the naturally associated 
super-algebra map $A \to T$, factors through $S$
by some $A\to S$ or 
by the naturally induced $\overline{A}\to S$. 

(2)
We claim that the super-ideal $I_A\, (=(A_1))$ of $A$ is nilpotent. 
Indeed, since $H$ is finitely generated by assumption, $\A\ot_BA\, (\simeq \A\ot H)$ is finitely
generated over $\A$. Since $B \to A^{\natural}$ is faithfully flat,
$A$ is finitely generated over $B$; this implies
that the super-ideal
$I_A$, being generated by finitely many odd elements, is nilpotent.

The claim above, combined with Part 1, shows that the $B$-super-algebra surjection $A \to \overline{A}$ 
(or in fact, $A_0\to \overline{A}$) splits.
The faithful flatness of $A$ over $B$ implies that
$\overline{A}$ is faithful flat over $B$.
\epf

Apply $\overline{A}\otimes_A$ to the isomorphism $\beta$ in \eqref{EBET} to obtain
\begin{equation}\label{EBETAA}
\overline{A}\otimes_BA\simeq \overline{A}\otimes H. 
\end{equation}
Then one sees that $A$ is a twisted form of the $B$-algebra $B\ot H$ in $\M^H$,
which is split by $\overline{A}$; this $\overline{A}$ is 
faithfully flat over $B$ by Lemma \ref{LBF} (2). 
Such twisted forms are classified by the 1st Amitsur cohomology set
\begin{equation}\label{EACH1}
H^1(\overline{A}/B,\ \G)
\end{equation}
for the affine super-group $\G$ represented by $H$.
This cohomology set is in fact constructed in the super situation; the construction
is analogous to the one, as given in \cite[Part V]{W}, of 
the ordinary cohomology set. One sees that
super-algebra maps 
$H \to \bigotimes_B^n \overline{A}$, 
$1 \le n \le 3$,
uniquely factor through $\overline{H}$.
Moreover, the complex for computing \eqref{EACH1} is naturally identified 
with the one for computing the ordinary 1st Amitsur cohomology set
\[
H^1(\overline{A}/B,\ \G_{\mathsf{ev}})
\]
for the affine group $\G_{\mathsf{ev}}$ represented by $\overline{H}$.
We see easily that the identification 
\[
H^1(\overline{A}/B,\ \G_{\mathsf{ev}})
\overset{\simeq}{\longrightarrow}
H^1(\overline{A}/B,\ \G)
\]
is realized by
\[ D \mapsto D \, \square_{\overline{H}}H, \]
where $D$ is a twisted form of $B\ot \overline{H}$ which is split by $\overline{A}$.
In particular, $A$ is (uniquely) of the form $D \, \square_{\overline{H}}H$. 
To be more precise, $A$ is isomorphic to $D \, \square_{\overline{H}}H$ through 
an isomorphism such as in \eqref{EIOT}, which is seen by \eqref{ETPDop}
to induce $\overline{A}\simeq D$.
Being a twisted form of $B\ot \overline{H}$, $\overline{A}\, (=D)$ is $\overline{H}$-Galois over $B$. 
This completes the proof of Proposition \ref{PPEB}.

\subsection{Step~3}\label{STHA3}
Under the situation at the beginning of Section \ref{STHA}, 
we assume that $H$ is finitely generated. In addition to (i), we assume
\begin{itemize}
\item[(ii)] the natural algebra surjection $B_0\to \overline{B}$ splits, and 
\item[(iii)] the super-ideal $(B_1)$ of $B$ is nilpotent,
\end{itemize}
as in Theorem \ref{THA}. We are going to complete the proof of the theorem.   

Let 
\[ F:= A/B_1A. \]
This is seen to be an algebra in $\M^H$. 
Given a purely even $\overline{B}$-super-module $N$, 
one sees
\[
N\ot_BA=(N\otimes_B\overline{B})\otimes_B A=N\ot_{\overline{B}}F. 
\]
It follows that $F$ is $\overline{B}$-faithfully flat, since $A$ is $B$-faithfully flat.
 
By applying $\overline{B}\, \otimes_B$ to the isomorphism $\beta$ in \eqref{EBET}, we obtain
\[ F \otimes_{\overline{B}}F \simeq F \otimes H. \]
Therefore, $F$ is $H$-Galois over the purely even $\overline{B}$. 
By Proposition \ref{PHS}, $F$ is $H$-smooth over $\overline{B}$. 
By (ii), the natural surjection $A \to F$ may be regarded as a $\overline{B}$-algebra
morphism in $\M^H$, which splits since its kernel is nilpotent by (iii). 
Choose arbitrarily a section $F \to A$, through which we can regard $A$ as an $F$-algebra
in $\M^H$, and in particular, as an object in $\M_F^H$. 

Notice that $\overline{F}=\overline{A}$. Then it follows by Proposition \ref{PPEB} that
\[
F \simeq \overline{A}\, \square_{\overline{H}}H,
\]
and $\overline{A}/\overline{B}$ is $\overline{H}$-Galois. 
The category equivalence $\M_F^H\approx \M_{\overline{B}}$ as
in \eqref{EPHI} shows
\[
A \simeq 
B\otimes_{\overline{B}}F\simeq B\otimes_{\overline{B}}(\overline{A}\square_{\overline{H}}H)
=(B\otimes_{\overline{B}}\overline{A})\square_{\overline{H}}H,
\]
as desired; see Remark \ref{RSmooth} (3) for the last canonical isomorphism.


\section{Proof of Theorem \ref{THB}}\label{STHB}

Let $R$ be a super-algebra.
Recall from \eqref{EBAR} that
we write $I_R=(R_1)$ and $\overline{R}=R/I_R$.
Define
\[ \operatorname{gr}R:=\bigoplus_{n\ge 0}I_R^n/I_R^{n+1}. \]
Thus we have
\[
(\operatorname{gr}R)(0)=\overline{R},\quad (\operatorname{gr}R)(1)=I_R/I_R^2,
\]
in particular. 
This $\operatorname{gr}R$ is a \emph{graded super-algebra}, by which we mean
an $\mathbb{N}$-graded algebra which, regarded to be $\mathbb{Z}_2$-graded 
by mod 2 reduction, is super-commutative. In other words, a \emph{graded super-algebra}
is a commutative algebra in the tensor category
\begin{equation}\label{EGM}
\grM_{k}=(\grM_{k}, \otimes, k)
\end{equation}
of the $\mathbb{N}$-graded vector spaces over $k$, equipped with the super-symmetry for
which objects are regarded to be $\mathbb{Z}_2$-graded by mod 2 reduction; for example,
exterior algebras, which have appeared so far, are such algebras.
We see
\[ \overline{\operatorname{gr}R}=\overline{R}. \]
A \emph{graded Hopf super-algebra} is a commutative Hopf algebra in the symmetric category
above. 

To prove Theorem \ref{THB} (in 3 steps), suppose that we are in the situation of the theorem.
Thus, $H$ is a finitely generated Hopf super-algebra, $A=(A,\rho)$ is a non-zero algebra
in $\M^H$, and $B=A^{\mathrm{co}\hspace{0.2mm}H}$.

\subsection{Step 1}\label{STHB1}
From the four assumptions (a)--(d) in Theorem \ref{THB}, we here assume the following three:
\begin{itemize}
\item[(b)] 
$\overline{A}/\overline{A}^{\mathrm{co}\hspace{0.2mm}\overline{H}}$ is $\overline{H}$-Galois;
\item[(c)] 
$\overline{H}$ is smooth as an algebra;
\item[(d)] 
$A$ is Noetherian and smooth.
\end{itemize}

\begin{lemma}\label{LHSA}
We have the following.
\begin{itemize}
\item[(1)] $\overline{A}^{\mathrm{co}\hspace{0.2mm}\overline{H}}$ is Noetherian and smooth.
\item[(2)] $\overline{A}$ is $\overline{H}$-smooth over $k$.
\end{itemize}
\end{lemma}
\pf
(1)\ By (d), $\overline{A}$ is Noetherian and smooth; 
for the smoothness, see the proof of Lemma \ref{LBF} (1).

In general, given a faithfully flat homomorphism $R \to S$ of algebras
over a field, we have that if $S$ is Noetherian and smooth, then $R$ is, as well. 
Indeed, $R$ is then clearly Noetherian. To see that it is smooth or equivalently,
geometrically regular, it suffices by base extension and localization to prove
that $R$ is regular, assuming that $R \to S$ is a flat local homomorphism  of
Noetherian local algebras with $S$ being regular. The desired regularity of $R$
follows by using
Serre's characterization of regularity as having finite global dimension.

Part 1 follows from this general fact applied to 
$\overline{A}^{\mathrm{co}\hspace{0.2mm}\overline{H}}\hookrightarrow \overline{A}$, which is faithfully flat
by (b).

(2)\
Everything here is purely even. We wish to prove that
an extension
\[ 0\to M \to T \to \overline{A} \to 0 \]
of commutative algebras 
in $\M^{\overline{H}}$ necessarily splits. This splits as a short 
exact sequence in $\M^{\overline{H}}$, since $M$ is in $\M^{\overline{H}}_{\overline{A}}$, and
is, therefore, injective in $\M^{\overline{H}}$ by (b); 
see the second paragraph of the proof of Lemma \ref{LSP}. 
Therefore, the algebra map 
$T^{\mathrm{co}\hspace{0.2mm}\overline{H}} \to 
\overline{A}^{\mathrm{co}\hspace{0.2mm}\overline{H}}$ restricted
to the $\overline{H}$-co-invariants is surjective, whence
it splits by Part 1 above. The desired result follows, since $\overline{A}$ is
$\overline{H}$-smooth over 
$\overline{A}^{\mathrm{co}\hspace{0.2mm}\overline{H}}$, as is seen from
Proposition \ref{PHS} applied to $\overline{A}/\overline{A}^{\mathrm{co}\hspace{0.2mm}\overline{H}}$; the last application is possible by (b) and (c). 
\epf

\begin{prop}\label{PBSM}
$A^{\mathrm{co}\hspace{0.2mm}\overline{H}}$ is a Noetherian smooth
super-algebra such that
\[
\overline{A^{\mathrm{co}\hspace{0.2mm}\overline{H}}}=\overline{A}^{\mathrm{co}\hspace{0.2mm}\overline{H}}.
\]  
\end{prop}
\pf
By (d), $I_A$ is nilpotent. By Lemma \ref{LHSA} (2) 
we have a section of the algebra projection $A \to \overline{A}$
in $\M^{\overline{H}}$, through which we can and we do regard $A$ as 
an $\overline{A}$-algebra in 
$\M^{\overline{H}}$. 
Let 
\begin{equation}\label{EPA}
P_A:=I_A/I_A^2\, (=(\operatorname{gr} A)(1)).
\end{equation}
Then this is an object in $\M_{\overline{A}}^{\overline{H}}$,
which is finitely generated projective $\overline{A}$-super-module by
Proposition \ref{PSMOO}. 
The natural projection
$I_A \to P_A$
is now a morphism in $\M_{\overline{A}}^{\overline{H}}$, and has a section 
essentially by the
Maschke-type Theorem \cite[Theorem 1]{D}; see the proof of Lemma \ref{LSPR}. 
The argument of Remark \ref{RSMOO}, slightly modified, shows that we have an isomorphism
\begin{equation}\label{EXI}
\xi:\wedge_{\overline{A}}(P_A)\overset{\simeq}{\longrightarrow} A 
\end{equation}
of $\overline{A}$-algebras in $\M^{\overline{H}}$.
By the category equivalence
$\M_{\overline{A}}^{\overline{H}}\approx \M_{\overline{A}^{\mathrm{co}\hspace{0.2mm}\overline{H}}}$
ensured by the assumption (b), we see that $P_A$ is of the form
\[ P_A=Q_A\ot_{\overline{A}^{\mathrm{co}\hspace{0.2mm}\overline{H}}} \overline{A},\ \,
\text{where}\ \,
Q_A=P_A^{\ co\overline{H}}. \]
This $Q_A$ is a finitely generated projective $\overline{A}^{\mathrm{co}\hspace{0.2mm}\overline{H}}$-module
since $\overline{A}$ is faithfully flat over $\overline{A}^{\mathrm{co}\hspace{0.2mm}\overline{H}}$. 
We thus have an isomorphism 
\begin{equation}\label{EWQA}
\wedge_{\overline{A}^{\mathrm{co}\hspace{0.2mm}\overline{H}}}(Q_A)
\ot_{\overline{A}^{\mathrm{co}\hspace{0.2mm}\overline{H}}}\overline{A}\simeq A
\end{equation}
of $\overline{A}$-algebras in $\M^{\overline{H}}$, which restricts to an isomorphism 
\begin{equation*}\label{EAcoHbar}
\wedge_{\overline{A}^{\mathrm{co}\hspace{0.2mm}\overline{H}}}(Q_A)\simeq A^{\mathrm{co}\hspace{0.2mm}\overline{H}}
\end{equation*}
of super-algebras over $\overline{A}^{\mathrm{co}\hspace{0.2mm}\overline{H}}$. This implies the desired result
by Proposition \ref{PSMOO}. 
\epf


\subsection{Step 2}\label{STHB2}
In what follows we assume, in addition,
\begin{itemize}
\item[(a)] 
the super-algebra map $A\otimes A\to A\otimes H$,
$a\otimes b \mapsto a\rho(b)$ is surjective.
\end{itemize}

Applying $\operatorname{gr}$ to everything we write so as
\[
\cA=\op{gr} A,\quad \cH=\op{gr} H. 
\]
Then $\cH$ is a finitely generated graded Hopf super-algebra. Let $\grM^{\cH}$
denote the symmetric category of the right $\cH$-comodules in $\grM_{k}$; see
\eqref{EGM}. Then $\cA=(\cA,\operatorname{gr}(\rho))$ is an algebra in $\grM^{\cH}$. In particular,
$\cH$ is a Hopf super-algebra, and $\cA$ is an algebra in $\M^{\cH}$. 

\begin{lemma}
These $\cH$ and $\cA$ satisfy those conditions which corresponds to
(a)--(d) for $H$ and $A$;
explicitly, they are the same (b), (c) and the following modified two:
\begin{itemize}
\item[(a$'$)] the super-algebra map $\cA\otimes \cA\to \cA\otimes \cH$,
$a\otimes b \mapsto a\operatorname{gr}(\rho)(b)$ is surjective.
\item[(d$'$)] $\cA$ is Noetherian and smooth.
\end{itemize}
\end{lemma}
\begin{proof}
By Proposition \ref{PSMOO}, (d$'$) follows from (d). 

Since the super-algebra map in (a$'$) is induced from the one in (a) with $\mathrm{gr}$ applied,
one sees that (a) implies (a$'$), in view of the following general fact. Given a surjection of super-algebras, 
$R \to S$, the induced map 
$\operatorname{gr} R \to \operatorname{gr} S$ of graded super-algebras is
surjective.
Indeed, one sees that $R/I_R^2$ (with $I_R=(R_1)$ as before) 
has the even component $(R/I_R^2)_0$
which is canonically isomorphic to $\overline{R}=R/I_R$, and has
$(R/I_R^2)_1=I_R/I_R^2$ as the odd component. It follows 
that the given super-algebra surjection induces surjections
\[
\overline{R}\to \overline{S},\quad I_R/I_R^2\to I_S/I_S^2,
\]
which are identified with the neutral and the 1st components of 
the induced, graded super-algebra map. The map is surjective since the graded super-algebras
are generated by the described components.
\end{proof}

Thus we may assume (a$'$), (b), (c) and (d$'$).
We now claim that
the conclusion of Theorem \ref{THB} holds for $\cH$ and $\cA$, as follows.

\begin{lemma}\label{LHB}
$\cA/\cA^{\mathrm{co}\hspace{0.2mm}\cH}$ is $\cH$-Galois, and $\cA^{\mathrm{co}\hspace{0.2mm}\cH}$ is Noetherian, 
smooth and such that
\begin{equation}\label{EcAcocH}
\overline{\cA^{\mathrm{co}\hspace{0.2mm}\cH}}=\overline{A}^{\mathrm{co}\hspace{0.2mm}\overline{H}}
\end{equation}
\end{lemma}
\pf
Let $W:=(H^+/(H^+)^2)_1$ as in \eqref{EW}. 
Then $\cH$ is
(canonically) isomorphic to the smash co-product $\overline{H}~\rcosmash\, \wedge(W)$
just as in \eqref{EJJW}; see \cite[(3.14)]{MT}, for example. 
In other words, we have the split short exact sequence
\[ \overline{H}\rightarrowtail \cH\twoheadrightarrow \wedge(W) \]
of (graded) Hopf super-algebras, where $\wedge(W)$ is supposed to contain 
all elements of $W$ as primitives. 
Note that $\cA$ turns into an algebra in $\M^{\wedge(W)}$
along the last surjection $\cH\twoheadrightarrow \wedge(W)$. 
The assumption (a$'$)
implies that the alpha map (see \eqref{EALP}) for this 
$\cA$ in $\M^{\wedge(W)}$ is surjective. 
Remark \ref{RHB} (for $H$), now applied to the co-Frobenius $\wedge(W)$, shows that 
$\cA/\cA^{\mathrm{co}\hspace{0.2mm}(\wedge(W))}$ is $\wedge(W)$-Galois, whence 
$\cA^{\mathrm{co}\hspace{0.2mm}(\wedge(W))} \to \cA$ 
is faithfully flat, in particular.

Applying Proposition \ref{PHA} (for $H$) to $\wedge(W)$,
one can choose a morphism
$\wedge(W)\to \cA$ of algebras in $\M^{\wedge(W)}$. 
With $\operatorname{gr}$ applied, this morphism may be supposed to be graded. 
It uniquely extends to an isomorphism  
\[
\cA^{\mathrm{co}(\wedge(W))} \otimes \wedge(W)\simeq \cA
\]
of $\cA^{\mathrm{co}(\wedge(W))}$-algebras in $\grM^{\wedge(W)}$, as is seen by using the category
equivalence
\begin{equation}\label{EEQQ}
\M^{\wedge(W)}_{\cA}\overset{\approx}{\longrightarrow}
\M_{\cA^{\mathrm{co}(\wedge(W))}},\quad M \mapsto M^{\mathrm{co}(\wedge(W))}
\end{equation}
ensured by
$\cA/\cA^{\mathrm{co}(\wedge(W))}$ being $\wedge(W)$-Galois. 
It follows that
\begin{equation}\label{EAbar}
\overline{\cA^{\mathrm{co}(\wedge(W))}}=\overline{A}. 
\end{equation}
Moreover, by (d$'$), $\cA^{\mathrm{co}(\wedge(W))}$ is Noetherian and smooth.

Notice that $\cA^{\mathrm{co}(\wedge(W))}$ is an algebra in $\grM^{\overline{H}}$, such that
\[
(\cA^{\mathrm{co}(\wedge(W))})^{\mathrm{co}\hspace{0.2mm} \overline{H}}
=\cA^{\mathrm{co}\hspace{0.2mm}\cH}.
\]
Apply the category equivalence \eqref{EEQQ} to the surjection in (a$'$). 
The result is 
the faithfully-flat base extension
$\cA\ot_{\cA^{\mathrm{co}(\wedge(W))}}$
of the alpha map 
\[
\cA^{\mathrm{co}(\wedge(W))}\ot \cA^{\mathrm{co}(\wedge(W))}\to \cA^{\mathrm{co}(\wedge(W))}\ot \overline{H}
\]
for $\cA^{\mathrm{co}(\wedge(W))}$ in $\grM^{\overline{H}}$.
Therefore, this alpha map is surjective.
By \eqref{EAbar} we can apply (the argument of proving) Proposition \ref{PBSM}
to $\cA^{\mathrm{co}(\wedge(W))}$. 
It results that $\cA^{\mathrm{co}\hspace{0.2mm}\cH}\, 
(=(\cA^{\mathrm{co}(\wedge(W))})^{\mathrm{co}\hspace{0.2mm}\overline{H}})$ is
Noetherian and smooth, and we have an isomorphism 
\begin{equation}\label{ENA}
\wedge_{\overline{A}^{\mathrm{co}\hspace{0.2mm}\overline{H}}}(N_A)
\otimes_{\overline{A}^{\mathrm{co}\hspace{0.2mm}\overline{H}}}
\overline{A}\simeq \cA^{\mathrm{co}(\wedge(W))}
\end{equation}
of $\overline{A}$-algebras in $\grM^{\overline{H}}$.
This isomorphism is analogous to \eqref{EWQA}, but it is now canonical.
Here $N_A$ denotes the $\overline{H}$-co-invariants in the 1st component 
of $\cA^{\mathrm{co}(\wedge(W))}$, or in other words, 
\begin{equation}\label{ENAPA}
N_A=P_A \cap \cA^{\mathrm{co}\hspace{0.2mm}\cH};
\end{equation}
see \eqref{EPA}. This $N_A$ is a finitely generated projective
$\overline{A}^{\mathrm{co}\hspace{0.2mm}\overline{H}}$-module. 
The isomorphism \eqref{ENA} shows \eqref{EcAcocH}.
It also shows that
$\cA^{\mathrm{co}(\wedge(W))}$ 
is injective in $\M^{\overline{H}}$, whence by Theorem \ref{TOBER} (1), 
$\cA^{\mathrm{co}(\wedge(W))}/\cA^{\mathrm{co}\hspace{0.2mm}\cH}$ is $\overline{H}$-Galois.
Therefore, $\cA^{\mathrm{co}\hspace{0.2mm}\cH}\to \cA$, being a composite of two faithfully flat homomorphisms,
is faithfully flat.

It remains to show that the beta map 
$\cA\otimes_{\cA^{\mathrm{co}\hspace{0.2mm}\cH}}\cA\to \cA\otimes \cH$ (see \eqref{EBET})
for $\cA$ in $\grM^{\cH}$ is bijective. After applying \eqref{EEQQ}
one has only to prove
that 
\[
\cA\otimes_{\cA^{\mathrm{co}\hspace{0.2mm}\cH}}\cA^{\mathrm{co}(\wedge(W))}\to 
\cA\otimes \overline{H}
\]
is bijective. 
But this is a base extension of the beta map 
\[
\cA^{\mathrm{co}(\wedge(W))}\otimes_{\cA^{\mathrm{co}\hspace{0.2mm}\cH}}
\cA^{\mathrm{co}(\wedge(W))}
\overset{\simeq}{\longrightarrow}
\cA^{\mathrm{co}(\wedge(W))}\otimes \overline{H}
\]
for $\cA^{\mathrm{co}(\wedge(W))}$ in $\grM^{\overline{H}}$, which is bijective since 
$\cA^{\mathrm{co}(\wedge(W))}/\cA^{\mathrm{co}\hspace{0.2mm}\cH}$ is $\overline{H}$-Galois. This completes the proof. 
\epf


\subsection{Step 3}\label{SHB3}

We have seen that $\cA/\cA^{\mathrm{co}\hspace{0.2mm}\cH}$ is an $\cH$-Galois extension
which satisfies the assumptions of Theorem \ref{THA}. 
Let $\cR:={}^{\mathrm{co}\hspace{0.2mm}\overline{H}}\cH$; this $\cR$ is seen to be an algebra in $\grM^{\cH}$. 
The theorem tells us that $\cA/\cA^{\mathrm{co}\hspace{0.2mm}\cH}$ arises (uniquely) from
an $\overline{H}$-Galois extension over $\cA^{\mathrm{co}\hspace{0.2mm}\cH}$. By Proposition \ref{PARI} we have
an algebra morphism $f : \cR \to \cA$ in $\M^{\cH}$, which can be chosen,
with $\operatorname{gr}$ applied, so as to be in $\grM^{\cH}$. 
Define
$D_f:= \cA/(f(\cR^+))$, just as in \eqref{EDF}, 
and let
$p : \cA \to D_f$
denote the natural projection; these are in $\grM^{\overline{H}}$. Let 
\[
\iota : \cA\overset{\simeq}{\longrightarrow}
D_f \, \square_{\overline{H}} \cH, 
\]
be the isomorphism of algebras in $\grM^{\cH}$, such as given by \eqref{EIOT}. 
Recall that this $\iota$ is induced from the composite
\[
\cA\to \cA\otimes \cH
\overset{p \otimes\operatorname{id}_{\cH}}{\longrightarrow} D_f \ot \cH, 
\]
where the first arrow indicates the $\cH$-co-action on $\cA$.
Using the isomorphism $\xi$ in \eqref{EXI}, 
define $q : A \to D_f$ to be the composite
\begin{equation}\label{EqADf}
q :A\overset{\xi^{-1}}{\longrightarrow} \wedge_{\overline{A}}(P_A)
\overset{\operatorname{gr}(\xi)}{\longrightarrow}\cA \overset{p}{\longrightarrow}D_f. 
\end{equation}
This is an algebra morphism in $\M^{\overline{H}}$ such that $\operatorname{gr}(q)=p$.
Note that the composite
\[
A\to A\otimes H
\overset{q \otimes\operatorname{id}_{H}}{\longrightarrow} D_f \ot H, 
\]
where the first arrow indicates the $H$-co-action on $A$, induces an algebra morphism
in $\M^H$,
\[
\eta : A \to D_f \, \square_{\overline{H}} H,
\]
which is indeed an isomorphism, 
since we see $\operatorname{gr}(\eta)=\iota$.

One sees from the isomorphism $\iota$ that $D_f$ is isomorphic, as an algebra in 
$\M^{\overline{H}}$, to 
$\cA^{\mathrm{co}(\wedge{W})}$. Therefore, $D_f$ is injective in $\M^{\overline{H}}$, as is seen again from
\eqref{ENA}. 
The isomorphism $\eta$ tells us that $A$ is injective in $\M^H$. 
By Theorem \ref{TOBER} (1), 
the assumption (a) implies that $A/A^{\mathrm{co}\hspace{0.2mm}H}$ is $H$-Galois. 
Moreover, $\eta$ and $\iota$ induce the isomorphisms
\begin{equation}\label{EADA}
A^{\mathrm{co}\hspace{0.2mm}H}\simeq (D_f)^{\mathrm{co}\hspace{0.2mm}\overline{H}}\simeq \cA^{\mathrm{co}\hspace{0.2mm}\cH}
\end{equation}
of super-algebras. More precisely, the composite 
\[
A\overset{\xi^{-1}}{\longrightarrow} \wedge_{\overline{A}}(P_A)
\overset{\operatorname{gr}(\xi)}{\longrightarrow}\cA 
\]
of  the first two isomorphisms in \eqref{EqADf} induces an isomorphism 
$A^{\mathrm{co}\hspace{0.2mm}H}\simeq \cA^{\mathrm{co}\hspace{0.2mm}\cH}$.
By Lemma \ref{LHB}, $A^{\mathrm{co}\hspace{0.2mm}H}$ is Noetherian, smooth and such that
\[
\overline{A^{\mathrm{co}\hspace{0.2mm}H}}=\overline{A}^{\mathrm{co}\hspace{0.2mm}\overline{H}}.
\]
This completes the proof of Theorem 
\ref{THB}. 

\begin{rem}\label{RBgrB}
By applying $\operatorname{gr}$ to \eqref{EADA} we obtain an isomorphism
\[
\operatorname{gr}(A^{\mathrm{co}\hspace{0.2mm}H})\ \simeq \ 
(\operatorname{gr}A)^{\mathrm{co}(\operatorname{gr}H)},
\]
which is indeed canonical, arising from the inclusion $A^{\mathrm{co}\hspace{0.2mm}H}\hookrightarrow A$ with $\operatorname{gr}$ applied. 
Let $B=A^{\mathrm{co}\hspace{0.2mm}H}$, as above. Then it follows from \eqref{ENA} that the isomorphism 
in degree 1 gives
\begin{equation}\label{EIBNA}
I_B/I_B^2= N_A.
\end{equation}
\end{rem}

\section{Proof of Theorem \ref{THMII}}\label{STHMII}
Suppose that we are in the situation of Theorem \ref{THMII}. We work in the geometrical 
setting, in which super-schemes are regarded as super-ringed spaces; see \cite[Sect.~4]{MZ1}. 
Let $\cO_{\X}$ denote the structure sheaf of $\X$, in particular. 
As a ringed space, the associated scheme 
$\X_{\mathsf{ev}}$ has the same underlying topological
space as $\X$, and the structure sheaf 
assigns the algebra $\overline{\cO_{\X}(Y)}$ to every open set $Y$.

\begin{lemma}\label{LTHBtoTHMII}
Under the additional assumption that $\X$ and 
$\X_{\mathsf{ev}}\tilde{/}\G_{\mathsf{ev}}$ are both affine, Theorem \ref{THMII} holds.
\end{lemma}
\pf 
Suppose that we are in the situation of the beginning of Section \ref{SSHG}. 
We now assume that the Hopf super-algebra $H$ is finitely generated and the 
algebra $A$ in $\M^H$ is Noetherian.
Note that the freeness of $\mathsf{G}$-action on $\X$
is equivalent to saying that the map $\alpha$ in \eqref{EALP} is an epimorphism of
super-algebras. Therefore,
the lemma follows from Theorem \ref{THB}
if we prove that $\alpha$ is necessarily surjective. But this follows from
the following general fact for $R\to S$ applied to $\alpha$. 
\epf

\begin{lemma}\label{LEPI}
An epimorphism $R \to S$ of super-algebras is necessarily surjective,
if the induced algebra map $\overline{R}\to \overline{S}$ is surjective and if 
$S$ is Noetherian. 
\end{lemma}
\pf
The assumptions imply that $S$, regarded as an $R$-super-algebra through the map,
is generated by finitely many nilpotent homogeneous elements in $I_R$, whence
it is finitely generated as an $R$-super-module. 
The desired surjectivity now follows by a familiar argument which uses the fact that
$R \to S$ is an epimorphism if and only if it, tensored $\ot_R\mathrm{id}_S$ with the identity map 
of $S$, turns into a bijection, $S=R \ot_R S \overset{\simeq}{\longrightarrow}S\ot_R S$; see
the proof of \cite[Proposition 2.6]{MSS}, for example. 
\epf

Let $\pi :\X_{\mathsf{ev}} \to \X_{\mathsf{ev}}\tilde{/}\G_{\mathsf{ev}}$ 
denote the natural morphism. 
Choose arbitrarily an affine open subset $U \ne \emptyset$ in 
$\X_{\mathsf{ev}}\tilde{/}\G_{\mathsf{ev}}$, which is necessarily Noetherian and smooth. 
Then, as is shown in 
\cite[Part I, 5.7, Page 83]{J}, $\pi^{-1}(U)$ is a $\G_{\mathsf{ev}}$-stable, open affine subscheme of $\X_{\mathsf{ev}}$ such that $\pi^{-1}(U)\tilde{/}\G_{\mathsf{ev}}=U$.
 
Recall that $\X$ and $\X_{\mathsf{ev}}$ have the same underlying topological space. 
One sees that $\pi^{-1}(U)$, regarded as an open subset of $\X$, is an open Noetherian 
super-subscheme
of $\X$, which we denote by $V$. 
The associated scheme is the $\pi^{-1}(U)$ in $\X_{\mathsf{ev}}$, which is affine,
as was seen above. Therefore, the $V$ in $\X$ is affine
by Zubkov's theorem \cite[Theorem 3.1]{Z}.
Clearly, the morphism of (affine) super-schemes
$V \times \G \to V \times V$ restricted from \eqref{EfALPHA} is a monomorphism, 
or in other words, the $\G$-action
on $V$ is free. By Lemma \ref{LTHBtoTHMII} we have 
$V \tilde{/}\G=\operatorname{Spec}(\cO_{\X}(V)^{\G})$, and this is Noetherian and smooth. Moreover, the original $U$ in 
$\X_{\mathsf{ev}}\tilde{/}\G_{\mathsf{ev}}$ coincides with 
$\operatorname{Spec} \overline{B}$, where we let $B:=\cO_{\X}(V)^{\G}$.

Let $Z$ denote the underlying topological space of the scheme $\X_{\mathsf{ev}}\tilde{/}\G_{\mathsf{ev}}$. 

\begin{lemma}
There uniquely exists a  sheaf $\cF$ of super-algebras on $Z$, such that 
\[
(U, \cF|_U)=\operatorname{Spec}(\cO_{\X}(V)^{\G})
\]
for every $U$ with $V =\pi^{-1}(U)\, (\subset \X)$, as above. 
\end{lemma}
\pf
Let $U\supset U'$ be non-empty affine open subsets of 
$\X_{\mathsf{ev}}\tilde{/}\G_{\mathsf{ev}}$. 
Define $V:=\pi^{-1}(U)$ and $V':=\pi^{-1}(U')$ in $\X$. The commutative diagram
\[
\begin{xy}
(0,0)   *++{V}  ="1",
(18,0)  *++{V'} ="2",
(0,-16)  *++{U} ="3",
(18,-16) *++{U'}="4",
{"2" \SelectTips{cm}{} \ar @{_{(}->} "1"},
{"4" \SelectTips{cm}{} \ar @{_{(}->}"3"},
{"1" \SelectTips{cm}{} \ar @{->}_{\pi|_{V}} "3"},
{"2" \SelectTips{cm}{} \ar @{->}^{\pi|_{V'}} "4"}
\end{xy}
\]
of affine super-schemes arises from the commutative diagram
\[
\begin{xy}
(0,0)   *++{A}  ="1",
(18,0)  *++{A'} ="2",
(0,-16)  *++{B} ="3",
(18,-16) *++{B'}="4",
{"1" \SelectTips{cm}{} \ar @{->}"2"},
{"3" \SelectTips{cm}{} \ar @{->}"4"},
{"3" \SelectTips{cm}{} \ar@{^{(}->} "1"},
{"4" \SelectTips{cm}{} \ar @{^{(}->} "2"}
\end{xy}
\]
of super-algebras, where we have set
\[ A =\cO_{\X}(V),\quad B=A^{\mathrm{co}\hspace{0.2mm}H}; \quad A'=\cO_{\X}(V'),\quad B'=A'^{\mathrm{co}\hspace{0.2mm}H}.\]
Choose arbitrarily a point $y \in U'$. There is a point $x \in V'$ such that $\pi(x)=y$.
Let 
\[
\mathfrak{P}^{(')}\in \operatorname{Spec}(A_0^{(')}),\quad
\mathfrak{p}^{(')}\in \operatorname{Spec}(B_0^{(')})
\]
be the primes which correspond to the point $x$ in $V^{(')}$ and to $y$ in $U^{(')}$,
respectively. To be precise, $\mathfrak{p}$ is first chosen from
$\operatorname{Spec}(\overline{B})$, and is then regarded to be in 
$\operatorname{Spec}(B_0)$ through the canonical identification of the two
sets of primes. This is the case for $\mathfrak{p}'$, too.  
We have the canonical
isomorphism
\begin{equation}\label{EAPBp1}
A_{\mathfrak{P}}\ = \ A'_{\mathfrak{P'}}
\end{equation}
between stalks of $\cO_{\X}$,  
and the analogous one 
\begin{equation}\label{EAPBp2}
\overline{B}_{\mathfrak{p}}\ = \ \overline{B'}_{\mathfrak{p'}}
\end{equation}
for the structure sheaf of
$\X_{\mathsf{ev}}\tilde{/}\G_{\mathsf{ev}}$. 
To prove the lemma it suffices to prove the canonical map 
\begin{equation*}\label{ECANO}
\mathrm{cano} : \ B_{\mathfrak{p}}\to B'_{\mathfrak{p}'}
\end{equation*}
is an isomorphism, which shows $\cO_{\operatorname{Spec}B}|_{U'}=
\cO_{\operatorname{Spec}B'}$. 
Notice that
the super-algebras $B_{\mathfrak{p}}$ and $B'_{\mathfrak{p}'}$ both
are Noetherian and smooth.

As is easily seen, localizations by the primes above are compatible
with the relevant constructions which will appear below; for example, we have 
$(\overline{B})_{\mathfrak{p}}=\overline{(B_{\mathfrak{p}})}$. We will use freely
such identifications. 

Just as in \eqref{EPA}, we write $P_R=I_R/I_R^2$ for $R= B^{(')}_{\mathfrak{p}^{(')}}$.
Let us return to the
isomorphism given in \eqref{ENA}, where we wrote $\cA=\operatorname{gr} A$. Localize 
the both sides at $\mathfrak{P}$, and restrict to the first components.
By using \eqref{EIBNA}, there results
\[
P_{B_{\mathfrak{p}}}\otimes_{\overline{B}_{\mathfrak{p}}}\overline{A}_{\mathfrak{P}}
=\cA^{\ co(\wedge(W))}_{\mathfrak{P}}(1). 
\]
Here we remark that since $\wedge(W)$ co-acts trivially on the $\overline{A}\, (=\cA(0))$
in $\cA$, $\cA^{\ co(\wedge(W))}_{\mathfrak{P}}$ makes sense, and 
its first component has appeared on the right-hand side. 
The analogous result for $A'$, $B'$ and 
$\cA':=\operatorname{gr} A'$,
combined with \eqref{EAPBp2} and $\cA_{\mathfrak{P}}=\cA'_{\mathfrak{P}'}$ which arises from \eqref{EAPBp1}, gives 
a canonical isomorphism
\begin{equation}\label{EPBPB}
P_{B_{\mathfrak{p}}}\ =\ P_{B'_{\mathfrak{p}'}}.
\end{equation}

Choose first a section $s$ of the super-algebra projection
$B_{\mathfrak{p}}\to \overline{B}_{\mathfrak{p}}$, and then a section $t$ 
of the projection
$I_{B_{\mathfrak{p}}}\to P_{B_{\mathfrak{p}}}$
in $\M_{\overline{B}_{\mathfrak{p}}}$,
where the $B_{\mathfrak{p}}$-super-module $I_{B_{\mathfrak{p}}}$ is regarded to be in
$\M_{\overline{B}_{\mathfrak{p}}}$ though the first chosen $s$. Define $s'$ and $t'$
to be the composites
\[
\overline{B'}_{\mathfrak{p}'}=\overline{B}_{\mathfrak{p}}\overset{s}{\longrightarrow}
B_{\mathfrak{p}}
\overset{\mathrm{cano}}{\longrightarrow} B'_{\mathfrak{p}'},\quad 
P_{B'_{\mathfrak{p}'}}=P_{B_{\mathfrak{p}}}\overset{t}{\longrightarrow} I_{B_{\mathfrak{p}}}\overset{\mathrm{cano}|_{I_{B_{\mathfrak{p}}}}}{\longrightarrow} I_{B'_{\mathfrak{p}'}},
\]
respectively. Clearly, these are sections of
$B'_{\mathfrak{p}'}\to \overline{B'}_{\mathfrak{p}'}$ and
of $I_{B'_{\mathfrak{p}'}}\to P_{B'_{\mathfrak{p}'}}$, respectively.
We have the following diagram:
\[
\begin{xy}
(0,0) *++{\wedge_{\overline{B}_{\mathfrak{p}}}(P_{B_{\mathfrak{p}}})}  ="1",
(26,0)  *++{B_{\mathfrak{p}}} ="2",
(0,-18)  *++{\wedge_{\overline{B'}_{\mathfrak{p}'}}(P_{B'_{\mathfrak{p}'}})} ="3",
(26,-18) *++{B'_{\mathfrak{p}'}}="4",
{"1" \SelectTips{cm}{} \ar @{->} "2"},
{"3" \SelectTips{cm}{} \ar @{->}"4"},
{"1" \SelectTips{cm}{} \ar @{->} "3"},
{"2" \SelectTips{cm}{} \ar @{->}^{\mathrm{cano}} "4"}
\end{xy}
\]
Here the horizontal arrow on the top (resp., on the bottom) indicates
the isomorphism (see \eqref{EWIR})
which arises from $s$ and $t$ (resp., from $s'$ and $t'$), and
the vertical arrow on the left-hand side indicates the isomorphism 
which arises from \eqref{EAPBp2} and \eqref{EPBPB}. 
It follows that the vertical $\mathrm{cano}$ on the right-hand side is an
isomorphism, as desired, since the diagram above is easily seen to be commutative.
\epf

We see that $\mathsf{Z}:=(Z,\cF)$ is a locally Noetherian, smooth super-scheme.
This represents $\X\tilde{/}\G$, and satisfies 
\[
\mathsf{Z}_{\mathsf{ev}}=\X_{\mathsf{ev}}\tilde{/}\G_{\mathsf{ev}}, 
\]
since it has locally these properties. This proves Theorem \ref{THMII}. 

\section*{Acknowledgments}
The first-named author was supported by JSPS~KAKENHI, Grant Numbers 17K05189 and 20K03552.

\end{document}